\documentclass[11pt, letterpaper]{elsarticle}
\usepackage{amsmath}
\usepackage{amssymb}
\usepackage{graphicx}
\usepackage{multirow}
\usepackage{latexsym}
\usepackage{epic}
\usepackage{url}
\usepackage{soul}
\usepackage{afterpage}
\usepackage[unicode, pdftex]{hyperref}

\usepackage{bbm}
\usepackage{setspace}
\usepackage{enumitem}
\usepackage{cleveref}
\usepackage[toc]{appendix}

\usepackage{tikz}
\usepackage{tikz,fullpage}
\usetikzlibrary{arrows,automata}
\usepackage{tkz-berge}
\usepackage{amsthm}
\usepackage[symbol]{footmisc}
\usepackage[mathscr]{euscript}
\usepackage[left=0.9in,top=0.9in,right=0.9in,bottom=0.9in,nohead]{geometry}
\DeclareMathOperator{\argmax}{argmax}
\DeclareMathOperator{\argmin}{argmin}

\makeatletter
\newtheorem*{rep@theorem}{\rep@title}
\newcommand{\newreptheorem}[2]{%
	\newenvironment{rep#1}[1]{%
		\def\rep@title{#2 \ref{##1}}%
		\begin{rep@theorem}}%
		{\end{rep@theorem}}}
\makeatother

\newreptheorem{example}{Example}
\newtheorem{theorem}{Theorem}

\newtheorem{proposition}{Proposition}
\newtheorem{definition}{Definition}

\newtheorem{example}{Example}

\usepackage{xcolor}
\usepackage[colorinlistoftodos,prependcaption,textsize=tiny]{todonotes}
\usepackage{xargs}

\usepackage{caption}
\usepackage{subcaption}

\usepackage{float}
\usepackage[section]{placeins}
\usepackage{arydshln}
\usepackage{latexsym}
\usepackage{epic}
\usepackage{amsmath,amsthm,amssymb}
\usepackage{graphicx}
\usepackage{url}
\usepackage{pgfplots}\usetikzlibrary{plotmarks}
\usepackage{soul}
\usepackage{afterpage}
\usepackage{cleveref}
\usepackage{bbm}
\usetikzlibrary{decorations.markings}
\usepackage[left=0.9in,top=0.9in,right=0.9in,bottom=0.9in,nohead]{geometry}
\usepackage{setspace}

\DeclareGraphicsExtensions{.pdf,.png,.jpg,.bmp}

\allowdisplaybreaks

\begin{document}
\begin{frontmatter}
\title{On a class of data-driven mixed-integer programming problems under uncertainty: a distributionally robust approach}

\author[label1]{Sergey S.~Ketkov\footnote[2]{Corresponding author. Email: sketkov@hse.ru; phone: +7 910 382-27-32.}}
\author[label1]{Andrei S.~Shilov}

\address[label1]{Laboratory of Algorithms and Technologies for Networks Analysis, National Research University \\Higher School of Economics (HSE), Bolshaya Pecherskaya st., 25/12, Nizhny Novgorod, 603155, Russia}

\begin{abstract}
In this study we analyze linear mixed-integer programming problems, in which the distribution of the cost vector is only observable through a finite training data set. In contrast to the related studies, we assume that the number of random observations for each component of the cost vector may vary. Then the goal is to find a \textit{prediction rule} that converts the data set into an estimate of the expected value of the objective function and a \textit{prescription rule} that provides an associated estimate of the optimal decision. We aim at finding the least conservative prediction and prescription rules, which satisfy some specified asymptotic guarantees as the sample size tends to infinity. We demonstrate that under some mild assumption the resulting vector optimization problems admit a Pareto optimal solution with some attractive theoretical properties. In particular, this solution can be obtained by solving a distributionally robust optimization~(DRO) problem with respect to all probability distributions with given component-wise relative entropy distances from the empirical marginal distributions. It turns out that the outlined DRO problem can be solved rather effectively whenever there exists an effective algorithm for the respective deterministic problem. In addition, we perform numerical experiments where the out-of-sample performance of the proposed approach is analyzed. 
\end{abstract}


\end{frontmatter}
\onehalfspace
\section{Introduction} \label{sec: intro}
 We consider a class of linear mixed-integer programming problems where the objective criterion is a function of both decision variables and uncertain problem parameters. The decision-maker cannot observe the nominal distribution of uncertain parameters, but has access to some finite training data set obtained from this distribution; see, e.g., the studies in \cite{Bertsimas2006, Esfahani2018}. 
Ideally, we endeavor to solve the following stochastic programming problem:
\begin{equation} \label{stochastic programming problem}
\min_{\mathbf{x} \in X} f(\mathbf{x}, \mathbb{Q}^*) := \mathbb{E}_{\mathbb{Q}^*} \{\gamma(\mathbf{c}, \mathbf{x})\},
\end{equation}
where $X \subseteq \mathbb{Z}^{n_1} \times \mathbb{R}^{n_2}$ is a linear mixed-integer set of feasible decisions, $\gamma(\mathbf{c}, \mathbf{x})$ is a given loss function and $\mathbf{c} \in \mathbb{R}^n$ with $n = n_1 + n_2$ is a vector of uncertain problem parameters (which is also referred to as a \textit{cost vector}) governed by some nominal probability distribution $\mathbb{Q}^*$. 
As outlined above, we assume that the nominal distribution $\mathbb{Q}^*$ is only observable through a finite data set. Furthermore, in contrast to the related literature, the data set is supposed to be \textit{incomplete}, i.e., the components of the cost vector $\mathbf{c}$ may be explored not to the same degree. The latter assumption is motivated by some online combinatorial optimization and machine learning problem settings discussed within Section \ref{subsec: approach and contribution}. 


From a computational perspective, if the nominal distribution $\mathbb{Q}^*$ is not available to the decision-maker, then we cannot even calculate the objective function value in (\ref{stochastic programming problem}) for a fixed decision~$\mathbf{x} \in X$. For this reason, the related solution approaches attempt to construct a fairly good approximation of the objective function value in (\ref{stochastic programming problem}) based on a set of available data.

\subsection{Related literature} \label{subsec: related literature}
In this section we briefly discuss how a finite set of random observations can be transformed to a solution of the stochastic programming problem (\ref{stochastic programming problem}). In the related literature this problem is usually addressed withing the framework of \textit{sample average approximation} (SAA) or \textit{distributionally robust optimization} (DRO). 

In the former approach the expected value of the loss function in (\ref{stochastic programming problem}) is approximated, for example, by the sample mean and the resulting function is optimized over the set $X$ of feasible decisions~\cite{Kleywegt2002}. In general, SAA methods provide deterministic formulations that enjoy strong \textit{asymptotic performance guarantees} due to the central limit theorem. In other words, the solutions obtained from SAA converge, in a sense, to the nominal solution of (\ref{stochastic programming problem}) as the sample size tends to infinity. It can also be argued that some robust modifications of sample average approximation methods may provide both asymptotic and finite sample performance guarantees; see, e.g., \cite{Bertsimas2018saa}.

The second approach to evaluate the expected value of the loss function in (\ref{stochastic programming problem}) is to construct an ambiguity set (or a family) of probability distributions, which is consistent with the available training data set. In other words, we must guarantee that the nominal distribution of uncertain parameters, $\mathbb{Q}^*$, belongs to the obtained ambiguity set with a sufficiently high probability. Then the decision-maker may solve a DRO problem, where the expected loss in (\ref{stochastic programming problem}) is optimized under the worst-case possible distribution of uncertain parameters; see, e.g., \cite{Delage2010, Goh2010, Wiesemann2014}.

In general, the problem of constructing an ambiguity set from data observations is considered by relatively many authors. For example, Delage and Ye \cite{Delage2010} design confidence sets for the support, mean and covariance matrix of uncertain parameters based on a finite set of independent training samples. Alternatively, several studies explore a distance metric in the space of probability distributions ``centered'' at the empirical distribution of the training samples; we refer to the studies in \cite{Bayraksan2015, Esfahani2018} and the references therein. Next, in our previous study \cite{Ketkov2021} we explore a situation where random observations of the cost vector $\mathbf{c}$ are not known \textit{exactly}, but reside within some specified affine manifolds. However, the modeling paradigm of \cite{Ketkov2021} is related to the shortest path problem and is restricted by a rather specific class of ambiguity sets. Eventually, we refer to Gupta \cite{Gupta2019}, where the author introduces a class of Bayesian ambiguity sets that are near-optimal in some specified asymptotic sense, i.e., when the sample size tends to infinity. 

The major drawback of the outlined solution techniques is that they decouple a \textit{prediction stage}, i.e., estimation of the objective function value in (\ref{stochastic programming problem}) for a fixed decision, and a \textit{prescription stage} where the obtained estimate is minimized over a given set of feasible decisions. In this regard, Van Parys~et~al.~\cite{VanParys2020} propose a \textit{meta-optimization problem}, which aims at finding, in a sense, optimal way of transforming the training data set to a solution of the stochastic programming problem (\ref{stochastic programming problem}). 

\looseness-1 The authors in \cite{VanParys2020} explore the case, in which the set of feasible decisions in (\ref{stochastic programming problem}) is continuous and the training data set is formed by independent observations of the \textit{vector} of uncertain~parameters. Van Parys~et~al.~\cite{VanParys2020} seek the \textit{least conservative} approximations of the expected loss in (\ref{stochastic programming problem}), for which some predefined asymptotic performance guarantees hold. It is proved that the optimal approximation in the aforementioned sense is unique and can be obtained by solving a DRO problem subject to all distributions within a ball with respect to the \textit{relative entropy distance}. Specifically, the radius of the ball controls the quality of asymptotic performance guarantees and the center of the ball is at the empirical distribution of the data. In addition, we refer to Sutter~et~al.~\cite{Sutter2020}, who applied the idea of Van~Parys~et~al.~\cite{VanParys2020} to non-i.i.d.~training data sets and more general stochastic programming~problems. 

\subsection{Our approach and contribution} \label{subsec: approach and contribution}
In this paper we propose another extension of the problem setting in \cite{VanParys2020} with respect to the construction of the training data set. The key modeling assumption of our study is that the number of observations for each particular component of the cost vector $\mathbf{c}$ is not necessarily the same and, thus, the overall data set is not complete. This assumption can be justified by a number of practical discrete optimization problem settings, in which the components of the cost vector $\mathbf{c}$ are explored not to the same degree.

At first, it is often the case that historical data is collected by trial and errors through multiple decision epochs. For example, in the stochastic multi-armed bandit problem \cite{Auer2002, Lai1985} the decision-maker has access to a finite set of actions with unknown expected rewards and attempts to maximize its cumulative reward over some finite time horizon. Naturally, in this problem setting it is more favorable for the decision-maker to explore the actions with higher expected rewards, which implies that some of the actions are explored more often than the others. The same idea can be applied to online combinatorial optimization; see, e.g., \cite{Audibert2014, Kveton2015}, where each action relates to a subset of ground items (for example, a routing decision in the context of network optimization problems corresponds to a sequence of arcs in a given network.) Then given a binary set of feasible decisions $X \subseteq \{0, 1\}^n$ the decision-maker may take a decision $\mathbf{x} \in X$ and observe only the costs associated with this decision, i.e., with nonzero components of $\mathbf{x}$.
	
Another possible application of incomplete data sets dates back to statistical analysis with missing data; see, e.g., \cite{Schafer2002} and the references therein. The existence of missing values in a training data set may depend on the data itself and hence, this fact is usually modeled as a probabilistic phenomenon. Depending on a concrete application missing values can be, for example, imputed based on some probabilistic model, ignored (parameters of the model are estimated based on the available data) or dropped from the analysis leaving only complete observations of the vector of uncertain parameters; see, e.g., \cite{Pigott2001}. With respect to discrete optimization problems, the source of missing values can be related to measurement errors, non-detects or simply a lack of information about some components of the cost vector $\mathbf{c}$.

In general, our modeling approach is motivated and similar to the approach of Van Parys~et~al.~\cite{VanParys2020}.
Formally, we seek an approximate solution of the stochastic programming problem~(\ref{stochastic programming problem}) by introducing \textit{prediction} and \textit{prescription rules}. A prediction rule converts the training data set into an estimate of the objective function value in (\ref{stochastic programming problem}) for some fixed decision $\mathbf{x} \in X$, whereas a prescription rule evaluates the optimal decision by minimizing the predicted value over the set $X$ of feasible decisions. The choice of prediction and prescription rules is induced by the following research~questions:
\begin{itemize}
	\item[\textbf{Q1}.] Are there prediction and prescription rules that can be applied to incomplete data sets and are also optimal in some predefined sense?
	\item[\textbf{Q2}.] If such optimal prediction and prescription rules exist, can we compute them rather effectively using off-the-shelf mixed-integer programming (MIP) solvers?
	\item[\textbf{Q3}.] How the sample size and the form of the nominal distribution affect the out-of-sample performance of the proposed model?
\end{itemize}

In order to address the first research question, \textbf{Q1}, we need to define some notion of optimality for the prediction and prescription rules. Following \cite{VanParys2020} the quality of these rules is estimated using an \textit{out-of-sample disappointment}, i.e., the probability that the nominal expected loss in (\ref{stochastic programming problem}), i.e., the expected loss under the nominal distribution $\mathbb{Q}^*$, is underestimated by the predicted expected loss. Next, we formulate \textit{prediction} and \textit{prescription problems}, which goal is to find the least conservative prediction and prescription rules, whose out-of-sample disappointment decays exponentially for every thinkable data-generating distribution $\mathbb{Q}^*$, as the sample size tends to infinity.

It turns out that the outlined prediction and prescription problems are, in fact, vector optimization problems over a functional space; see, e.g., \cite{Jahn2009}. For this reason, we distinguish between \textit{strongly} and \textit{weakly optimal} solutions. A strongly optimal solution of a vector optimization problem is strictly better than all other feasible solutions with respect to some specified partial order and a weakly optimal solution, in turn, is not dominated by any other feasible solution. We prove that under some additional assumptions on the objective criterion in (\ref{stochastic programming problem}) and the nominal distribution $\mathbb{Q}^*$,  there exists a rather broad class of weakly optimal solutions for both prediction and prescription problems. In particular, one of these solutions can be obtained by solving a DRO~problem subject to all probability distributions with given \textit{component-wise} relative entropy distances from the empirical \textit{marginal} distributions. We demonstrate that this solution inherits some properties of a strongly optimal solution and thoroughly analyze its theoretical properties.

With respect to the second research question, $\mathbf{Q2}$, our solution procedure for the aforementioned DRO problem is divided into two consecutive stages. In the first stage we retrieve the worst-case expected costs by solving univariate convex optimization problems for each component of the cost vector; in the second stage we solve a deterministic version of the underlying mixed-integer programming (MIP) problem. As a result, the DRO problem can be solved rather effectively whenever there exists an effective algorithm for the nominal MIP problem.

One natural limitation of our approach is that the objective criterion in (\ref{stochastic programming problem}) is assumed to be linear in the uncertain parameters. Intuitively, under this assumption we may shrink the probability space to the set of univariate distributions and consider component-wise ambiguity sets. It can also be argued that a straightforward implementation of the DRO approach from \cite{VanParys2020} with a discrete set of feasible decisions results in a substantially non-linear MIP problem, which cannot be solved at hand using off-the-shelf solvers; see Section \ref{subsec: benchmark approaches} for more details. 

Finally, we consider the research question \textbf{Q3} in our numerical study, where the out-of-sample performance of the proposed model in analyzed with respect to several classes of combinatorial optimization problems.
From the practical perspective, we demonstrate that solutions with a reasonably good quality can be obtained whenever the items with lower expected costs are observed more often than the items with higher expected costs. The latter observation seems to be relevant, e.g., if the historical data is collected by trial and errors through multiple decision epochs; recall our discussion of the online learning problem settings at the beginning of this section. In addition, we compare the DRO approach with several benchmark solution approaches based on measure concentration inequalities or truncation of the training data set and using the model of Van Parys~et~al.~\cite{VanParys2020}.

The remaining structure of the paper is summarized as follows:
\begin{itemize}
	\item In Section \ref{sec: problem} we introduce our prediction and prescription problems (that are motivated and similar to the related meta-optimization problems in \cite{VanParys2020}) under the assumption that the data set is incomplete. We also define problem-specific asymptotic and finite sample performance guarantees.
	\item Section \ref{sec: solution techniques} provides some analysis of the prediction and prescription problems. In particular, we focus on the quality and computational tractability of the proposed DRO approach.
	\item In Section \ref{sec: comp study} we provide numerical experiments with applications to data-driven shortest path and unweighted knapsack problems. The DRO approach is also compared with several benchmark approaches in terms of their average out-of-sample performance.
	\item Finally, Section~\ref{sec: conclusion} concludes the paper and outlines possible directions for future research.
\end{itemize}
\noindent \textbf{Notation.} All vectors and matrices are labelled by bold letters. 
The natural logarithm for some $q \in \mathbb{R}_{>0}$ is denoted as $\ln(q)$. Furthermore, we adapt the conventions $0 \ln(\frac{0}{q}) = 0$ and $q' \ln(\frac{q'}{0}) = \infty$ for any $q, q' \in \mathbb{R}_{>0}$.
The space of all probability distributions with a support $\mathcal{S} \subseteq \mathbb{R}^{n}$, $n \in \mathbb{Z}_{>0}$, is denoted by 
$\mathcal{Q}_0(\mathcal{S})$ and for any distribution~$\mathbb{Q} \in \mathcal{Q}_0(\mathcal{S})$ we denote by $\mathbf{Q} =~(\mathbb{Q}_1, \ldots, \mathbb{Q}_n)^\top$ a vector of the associated marginal distributions. We use $\mathbbm{1}\{Z\}$ as an indicator of some logical statement~$Z$. Finally, the probability of $Z$ with respect to some distribution $\mathbb{Q} \in \mathcal{Q}_0(\mathcal{S})$ is referred to as $\Pr_{\mathbb{Q}}\{Z\}$ (the index $\mathbb{Q}$ is omitted, if the data-generating distribution is clear from the context). 

\section{Problem formulation} \label{sec: problem}
\subsection{Modeling assumptions and terminology} \label{subsec: terminology}
As briefly outlined in Section \ref{sec: intro}, we consider the stochastic programming problem (\ref{stochastic programming problem}) with a mixed-integer set of feasible decisions $X \subseteq \mathbb{Z}^{n_1} \times \mathbb{R}^{n_2}$ and a random cost vector $\mathbf{c} \in \mathbb{R}^{n_1 + n_2}$. For simplicity of exposition we set $n := n_1 + n_2$ and define a set of indexes $\mathcal{A} := \{1, \ldots, n\}$. The cost vector $\mathbf{c}$ is assumed to be governed by some unknown nominal distribution $\mathbb{Q}^* \in \mathcal{Q}_0(\mathbb{R}^{n})$, which is only observable through a finite data set 
\begin{equation} \label{data set}
\widehat{\mathbf{C}} := \{(\widehat{c}_{a,1}, \ldots, \widehat{c}_{a,T_a})^\top, a \in \mathcal{A}\},
\end{equation}
Specifically, for each component $c_a$, $a \in \mathcal{A}$, of the cost vector $\mathbf{c}$ we presume existence of $T_a \in \mathbb{Z}_{>0}$ independent identically distributed (i.i.d.) samples from the related marginal distribution~$\mathbb{Q}^*_a$. 


Throughout the paper we make the following modeling assumptions:
\begin{itemize}
\item[\textbf{A1.}] Each component of the cost vector $\mathbf{c}$ is strictly positive and has a finite discrete support, i.e., 
\begin{equation} \label{eq: support}
c_a \in \mathcal{S}_a := \{z_{a,1}, \ldots, z_{a,d_a}\}
\end{equation} 
for some $z_{a,i} > 0$, $i \in \{1, \ldots, d_a\}$ and $a \in \mathcal{A}$. 
\item[\textbf{A2.}] Any feasible decision is nontrivial, nonnegative and bounded, i.e., $\mathbf{0} \notin X$ and for every $\mathbf{x} \in X$ there exists some $\mathbf{u} \in \mathbb{R}_+^n$ such that $x_a \in [0, u_a]$, $a \in \mathcal{A}$. 
\item[\textbf{A3.}] We assume that $\gamma(\mathbf{c}, \mathbf{x}) = \mathbf{c}^\top \mathbf{x}$ in the definition of the stochastic optimization problem (\ref{stochastic programming problem}).
\end{itemize}

The motivation behind Assumption \textbf{A1} is two-fold. First, it is argued in \cite{VanParys2020} that the case of continuous support requires more subtle mathematical techniques and makes the problem of finding optimal prediction and prescription rules substantially more difficult. For this reason, we leave an extension of our problem setting to the case of continuous support as a future research direction. Secondly, Assumption \textbf{A1} indicates that each component of the cost vector $\mathbf{c}$ has its individual support and, hence, there is \textit{no initial information} about some correlation between the components. The need of individual support sets can be motivated in the context of online combinatorial optimization \cite{Audibert2014, Kveton2015}, where the expected rewards are estimated individually for each action and there is no correlation information involved. 
On the other hand, the related assumption of Van Parys~et~al.~\cite{VanParys2020} requires complete knowledge of the support of the cost vector $\mathbf{c}$. 
Put differently, even in the absence of the training data set, the decision-maker in \cite{VanParys2020} \textit{has some initial information} about possible relations between the components of vector $\mathbf{c}$.

Next, Assumption \textbf{A2} is rather standard in the discrete optimization literature and holds for the majority of practical decision-making problems; see, e.g., \cite{Conforti2014}. Finally, Assumptions \textbf{A1} and \textbf{A3} imply that the objective function in (\ref{stochastic programming problem}) is completely defined by the marginal probability distributions $\mathbb{Q}^*_a$, $a \in \mathcal{A}$, induced by the nominal distribution $\mathbb{Q}^*$. In other words, we observe that
\begin{equation} \label{eq: objective function linearity}
f(\mathbf{x}, \mathbb{Q}^*) = \mathbb{E}_{\mathbb{Q}^*} \{\gamma(\mathbf{c}, \mathbf{x})\} = \sum_{a \in \mathcal{A}} \mathbb{E}_{\mathbb{Q}^*_a}\{c_a\} x_a
\end{equation}
 Furthermore, if the loss function $\gamma(\mathbf{c}, \mathbf{x})$ does not admit such a component-wise decomposition, then we cannot even compute the empirical mean of $\gamma(\mathbf{c}, \mathbf{x})$ as the associated training data set is incomplete.   
At the same time, some more complicated forms of the objective criteria in (\ref{stochastic programming problem}) are briefly discussed within Section \ref{sec: conclusion}. 


Below we provide some technical definitions that are used to formulate our prediction and prescription problems. Admittedly, our definitions are somewhat similar to the terminology of~Van~Parys~et~al. \cite{VanParys2020} but account incomplete knowledge of the training data set. 

First, in view of Assumption \textbf{A1}, the set of all marginal probability distributions for each $a \in \mathcal{A}$ can be defined as:
\begin{equation} \label{marginal ambiguity set}
\mathcal{Q}_a := \Big\{ \mathbb{Q} \in \mathcal{Q}_0(\mathcal{S}_a): \; 0 \leq q_{a,i} \leq 1 \quad \forall i \in \{1, \ldots, d_a\}, \; \sum_{i = 1}^{d_a} q_{a,i} = 1 \Big\},
\end{equation}
where $q_{a,i}$ denotes the probability that $c_a$ equals to $z_{a,i}$, $i \in \{1, \ldots, d_a\}$. Furthermore, taking into account the form of the objective function (\ref{eq: objective function linearity}) it suffices to consider the following ambiguity set of joint distributions: 
\begin{equation} \label{eq: set of all possible distributions}
\mathcal{Q} := \Big\{ \mathbb{Q} \in \mathcal{Q}_0(\mathcal{S}_1 \times \ldots \times \mathcal{S}_n): \mathbb{Q}_a \in \mathcal{Q}_a \quad \forall a \in \mathcal{A} \Big\}
\end{equation} 
As outlined earlier, for any joint distribution $\mathbb{Q} \in \mathcal{Q}$ we denote by 
$\mathbf{Q} = (\mathbb{Q}_{1}, \ldots, \mathbb{Q}_{n})^\top$ the vector of associated marginal distributions. 

The key idea of our theoretical analysis is to find, in a sense, optimal approximations of the nominal expected loss $f(\mathbf{x}, \mathbb{Q}^*)$ and an optimal~decision $\mathbf{x}^*(\mathbb{Q}^*) \in \argmin_{\mathbf{x} \in X} f(\mathbf{x}, \mathbb{Q}^*)$ in (\ref{stochastic programming problem}) by some functions of \textit{empirical marginal distributions}.

\begin{definition}[\textbf{Empirical marginal distributions}] \label{def: empirical marginal distribution} \upshape
	A univariate probability distribution $\widehat{\mathbb{Q}}_{a,T_a} \in \mathcal{Q}_a$ such that 
	\begin{equation} \label{eq: empirical marginal distributions}
	\widehat{q}_{a,i} := {\Pr}\{c_a = z_{a,i}\} = \frac{1}{T_a}\sum_{j = 1}^{T_a}\mathbbm{1}\{\widehat{c}_{a,j} = z_{a,i}\} \quad \forall i \in \{1, \ldots, d_a\}
	\end{equation}
	is referred to as an \textit{empirical marginal distribution} of $c_a$. \vspace{-10mm}\flushright$\square$
\end{definition}
\noindent In other words, by leveraging the available data set (\ref{data set}) we construct empirical estimators $\widehat{\mathbb{Q}}_{a,T_a}$ of the nominal marginal distributions $\mathbb{Q}^*_a$ for each $a \in \mathcal{A}$. Any vector of empirical marginal distributions that satisfies equation (\ref{eq: empirical marginal distributions}) is denoted as
$$\widehat{\mathbf{Q}}(T_1, \ldots, T_n) := (\widehat{\mathbb{Q}}_{1.T_1}, \ldots, \widehat{\mathbb{Q}}_{n,T_n})^\top,$$
where the dependence on $T_1, \ldots, T_n$ is usually omitted in order to streamline the notations.  

As a result, the definitions of prediction and prescription rules can be formalized as follows; we also refer to the related definition in \cite{VanParys2020}.  
\begin{definition} [\textbf{Prediction and prescription rules}] \upshape \label{def: data-driven predictor and prescriptor}
	 
	\looseness-1 A function $\hat{f}: X \times \mathcal{Q}_1 \times \ldots \times~\mathcal{Q}_n \rightarrow \mathbb{R}$ is called a data-driven \textit{prediction rule}, if 
	$\hat{f}(\mathbf{x}, \widehat{\mathbf{Q}})$ is used as an approximation of $f(\mathbf{x}, \mathbb{Q}^*)$. Furthermore, a function $\hat{x}: \mathcal{Q}_1 \times \ldots \times \mathcal{Q}_n \rightarrow X$ is an associated data-driven \textit{prescription rule}, if $\hat{x}(\mathbf{Q}) \in \argmin_{\mathbf{x} \in X} \hat{f}(\mathbf{x}, \mathbf{Q})$ for any feasible $\mathbb{Q} \in \mathcal{Q}$ and $\hat{x}(\widehat{\mathbf{Q}})$ is used as an approximation of $\mathbf{x}^*(\mathbb{Q}^*)$. \vspace{-8.75mm}\flushright$\square$
\end{definition}

Following \cite{VanParys2020} we assess the quality of approximation by using an \textit{out-of-sample disappointment}.

\begin{definition} [\textbf{Out-of-sample disappointment}] \upshape \label{def: out-of-sample prediction disappointment}
	For any prediction rule $\hat{f}$ and feasible decision $\mathbf{x} \in X$ an \textit{out-of-sample prediction disappointment} is given by:
	\begin{equation} \label{eq: out-of-sample prediction}
  {\Pr}_{\mathbb{Q}^*}\Big\{f(\mathbf{x}, \mathbb{Q}^*) > \hat{f}(\mathbf{x}, \widehat{\mathbf{Q}}) \Big \},
	\end{equation}
	 Similarly, for a given prediction-prescription pair $(\hat{f}, \hat{x})$
	\begin{equation} \label{eq: out-of-sample prescription}
	{\Pr}_{\mathbb{Q}^*}\Big\{f(\hat{x}(\widehat{\mathbf{Q}}), \mathbb{Q}^*) > \hat{f}(\hat{x}(\widehat{\mathbf{Q}}), \widehat{\mathbf{Q}}) \Big \}
	\end{equation}
	is reffered to as an \textit{out-of-sample prescription disappointment}. \vspace{-10mm}\flushright$\square$
\end{definition}
\looseness-1 The out-of-sample disappointment for fixed $\mathbf{x} \in X$ quantifies the probability that the nominal expected loss $f(\mathbf{x}, \mathbb{Q}^*)$ is underestimated by the predicted expected loss $\hat{f}(\mathbf{x}, \widehat{\mathbf{Q}})$. It can also be argued that, if a decision-maker aims at minimizing its loss, then underestimated losses are usually more harmful than overestimated losses; see, e.g., \cite{Sniedovich2011} for some basic concepts of worst-case analysis in classical decision theory.

As outlined in Section \ref{subsec: approach and contribution}, we focus on a class of prediction rules that cater an \textit{exponential decay rate} for the out-of-sample prediction disappointment irrespective of a decision $\mathbf{x} \in X$ and a data-generating distribution $\mathbb{Q}^* \in \mathcal{Q}$, as the sample size tends infinity. In other words, we set $T_{min} := \min_{a \in \mathcal{A}}T_a$ we exploit the following \textit{asymptotic} and \textit{finite sample performance guarantees}:
\begin{itemize}
	\item \looseness-1\textbf{Asymptotic guarantee.} The out-of-sample prediction disappointment (\ref{eq: out-of-sample prediction}) decays exponentially at some rate $r > 0$ when $T_{min}$ goes to infinity, i.e.,
	\begin{align} \label{eq: asymptotic guarantee} \tag{$\mathcal{AG}$}
	& \limsup_{T_{min} \rightarrow +\infty} \frac{1}{T_{min}} \ln \Big(\Pr\Big\{f(\mathbf{x}, \mathbb{Q}^*) > \hat{f}(\mathbf{x}, \widehat{\mathbf{Q}}(T_1, \ldots, T_n)) \Big \} \Big) \leq -r \quad \forall \mathbf{x} \in X, \; \forall \mathbb{Q}^* \in \mathcal{Q}
	\end{align}
	\item \textbf{Finite sample guarantee.}
	Let $g(T)$ be a given function that decays exponentially at some rate $r > 0$ as $T$ goes to infinity. Then for any $T_{a} \in \mathbb{Z}_{>0}$, $a \in \mathcal{A}$, the out-of-sample prediction disappointment is bounded from above by $g(T_{min})$, i.e.,
	\begin{align} \label{eq: finite sample guarantee} \tag{$\mathcal{FG}$}
	& \Pr\Big\{f(\mathbf{x}, \mathbb{Q}^*) > \hat{f}(\mathbf{x}, \widehat{\mathbf{Q}}(T_1, \ldots, T_n)) \Big \} \leq g(T_{min}) \quad \forall \mathbf{x} \in X, \; \forall \mathbb{Q}^* \in \mathcal{Q}, \; \forall T_{a} \in \mathbb{Z}_{>0}, \; a \in \mathcal{A}
	\end{align}
\end{itemize}

In the following example we demonstrate that the finite sample guarantee (\ref{eq: finite sample guarantee}) holds for a rather simple class of prediction rules. That is, we may approximate the expected costs in (\ref{stochastic programming problem}) by the sum of empirical mean and some positive constants.

\begin{example}[\textbf{Hoeffding bounds}] \upshape \label{example 1}
Let $\underline{z}_a = \min_{i \in \{1, \ldots, d_a\}} z_{a, i}$, $\overline{z}_a = \max_{i \in \{1, \ldots, d_a\}} z_{a, i}$ for each $a \in \mathcal{A}$ and $\mathcal{A}_\mathbf{x} := \{a \in \mathcal{A}: x_a >~0\} \neq \varnothing$; recall Assumption \textbf{A2}.
We assume that the prediction rule $\hat{f}$ for some fixed $\mathbf{x} \in X$ is defined~as:
\begin{equation} \label{eq: sample-average predictor + epsilon}
\hat{f}(\mathbf{x}, \widehat{\mathbf{Q}}) = \sum_{a \in \mathcal{A}} \min\Big\{\frac{1}{T_a}\sum_{j = 1}^{T_a} \widehat{c}_{a,j} + \varepsilon_a; \overline{z}_a \Big\} x_a 
\end{equation}	
where $\varepsilon_a > 0$ are some positive constants defined below and $\overline{z}_a$ serves as an upper bound on the expected costs for each $a \in \mathcal{A}$.  Using equation (\ref{eq: objective function linearity}) we observe that the out-of-sample prediction disappointment (\ref{eq: out-of-sample prediction}) can be bounded from above as: 
\begin{equation}
\begin{gathered} \allowdisplaybreaks \label{ex 1: proof}
\Pr\Big\{f(\mathbf{x}, \mathbb{Q}^*) > \hat{f}(\mathbf{x}, \widehat{\mathbf{Q}}) \Big \} = \Pr\Big\{\sum_{a \in \mathcal{A}} \Big(\mathbb{E}_{\mathbb{Q}^*_a}\{c_a\} - \min\Big\{\frac{1}{T_a}\sum_{j = 1}^{T_a} \widehat{c}_{a,j} + \varepsilon_a; \overline{z}_a \Big\}\Big) x_a > 0 \Big\} \leq \\
\Pr\Big\{\bigvee_{a \in \mathcal{A}_\mathbf{x}} \Big(\mathbb{E}_{\mathbb{Q}^*_a}\{c_a\} - \min\Big\{ \frac{1}{T_a}\sum_{j = 1}^{T_a} \widehat{c}_{a,j} + \varepsilon_a; \overline{z}_a \Big\} > 0 \Big) \Big\} \leq 
\\ 
\sum_{a \in \mathcal{A}_\mathbf{x}} \Pr\Big\{ \mathbb{E}_{\mathbb{Q}^*_a}\{c_a\} > \min\Big\{ \frac{1}{T_a}\sum_{j = 1}^{T_a} \widehat{c}_{a,j} + \varepsilon_a; \overline{z}_a \Big\}\Big\} \leq \\ \sum_{a \in \mathcal{A}_\mathbf{x}} \Pr\Big\{ \Big( \mathbb{E}_{\mathbb{Q}^*_a}\{c_a\} > \frac{1}{T_a}\sum_{j = 1}^{T_a} \widehat{c}_{a,j} + \varepsilon_a \Big) \bigvee \Big( \mathbb{E}_{\mathbb{Q}^*_a}\{c_a\} > \overline{z}_a \Big)\Big\} \leq 
\\
\sum_{a \in \mathcal{A}_\mathbf{x}}\Big( \Pr\Big\{ \mathbb{E}_{\mathbb{Q}^*_a}\{c_a\} > \frac{1}{T_a}\sum_{j = 1}^{T_a} \widehat{c}_{a,j} + \varepsilon_a \Big\} + \Pr\Big\{ \mathbb{E}_{\mathbb{Q}^*_a}\{c_a\} > \overline{z}_a \Big\}\Big) = 
\\
\sum_{a \in \mathcal{A}_\mathbf{x}} \Pr\Big\{ \mathbb{E}_{\mathbb{Q}^*_a}\{c_a\} > \frac{1}{T_a}\sum_{j = 1}^{T_a} \widehat{c}_{a,j} + \varepsilon_a \Big\}
\end{gathered}
\end{equation}

The first inequality in (\ref{ex 1: proof}) exploits Assumption \textbf{A2} and can be checked by contradiction, while the second inequality is implied by the standard union bound. In the following we use, respectively, the definition of minimum, the union bound and the fact that $\mathbb{E}_{\mathbb{Q}^*_a}\{c_a\} \leq \overline{z}_a$ for every $a \in \mathcal{A}$.

\looseness-1 Next, for any given decay rate $r > 0$ and $\delta_a \in (0, 1)$, $a \in \mathcal{A}$, such that $\sum_{a \in \mathcal{A}}\delta_a = 1$ we may set
\begin{equation} \label{eq: eps hoeffding}
\varepsilon_a = (\overline{z}_a - \underline{z}_a)\sqrt{\frac{rT_{min} - \ln \delta_a}{2T_a}}
\end{equation} 
By using Hoeffding's inequality \cite{Hoeffding1994} for the sum of bounded i.i.d. random variables we observe that
\begin{equation} \nonumber 
\Pr\Big\{ \mathbb{E}_{\mathbb{Q}^*_a}\{c_a\} - \frac{1}{T_a}\sum_{j = 1}^{T_a} \widehat{c}_{a,j} > - \varepsilon_a \Big\} \leq \exp\Big(-2 T_a \Big(\frac{ \varepsilon_a}{\overline{z}_a - \underline{z}_a}\Big)^2 \Big) = \delta_a e^{-rT_{min}}
\end{equation}
and, thus, for any $T_1, \ldots, T_n \in \mathbb{Z}_{>0}$
\begin{equation} \nonumber
\begin{gathered} \allowdisplaybreaks \label{ex 1: proof 2}
\Pr\Big\{f(\mathbf{x}, \mathbb{Q}^*) > \hat{f}(\mathbf{x}, \widehat{\mathbf{Q}}) \Big \} \leq
\sum_{a \in \mathcal{A}_\mathbf{x}} \delta_a e^{-rT_{min}} \leq e^{-rT_{min}} \end{gathered}
\end{equation}
The upper bounding function $g(T_{min}) := e^{-rT_{min}}$ has an exponential decay rate $r$, i.e., 
\begin{equation} \nonumber
\limsup_{T_{min} \to +\infty}\frac{1}{T_{min}} \ln g(T_{min}) = -r 
\end{equation}
We conclude that the prediction rule (\ref{eq: sample-average predictor + epsilon}) with the parameters $\varepsilon_a$, $a \in \mathcal{A}$, given by equation (\ref{eq: eps hoeffding}) satisfies both the finite-sample and asymptotic performance guarantees, (\ref{eq: finite sample guarantee}) and (\ref{eq: asymptotic guarantee}). \vspace{-7mm}\flushright$\square$
\end{example}

$ $\\
\subsection{Prediction and prescription problems} \label{subsec: prediction and prescription problems}

In this section we introduce prediction and prescription problems, which are also motivated by the related formulations in \cite{VanParys2020}.
First, we impose a partial order on the space of prediction rules. That is, a prediction rule $\hat{f}_1$ not strictly dominates a prediction rule $\hat{f}_2$, i.e., $\hat{f}_1 \preceq \hat{f}_2$, if and only if
\begin{equation} \label{eq: partial order prediction} \nonumber
\hat{f}_1(\mathbf{x}, \mathbb{Q}_1, \ldots, \mathbb{Q}_m) \leq \hat{f}_2(\mathbf{x}, \mathbb{Q}_1, \ldots, \mathbb{Q}_m) \quad \forall \mathbf{x} \in X, \; \forall \mathbb{Q} \in \mathcal{Q}
\end{equation}
Also, let $\mathcal{F}$ be a set of all real-valued functions on the set $X \times \mathcal{Q}_1 \times \ldots \times \mathcal{Q}_m$. Then a \textit{prediction problem} can be expressed as:
\allowdisplaybreaks\begin{subequations}
	\begin{align} \label{prediction problem} \tag{$\mathcal{P}$}
	& \min_{\hat{f}(\cdot) \in \mathcal{F}} \hat{f}(\cdot) \quad (\mbox{ w.r.t. }\preceq) \\
	\mbox{\upshape s.t.}
	& \limsup_{T_{min} \rightarrow +\infty} \frac{1}{T_{min}} \ln\Big(\Pr\Big\{f(\mathbf{x}, \mathbb{Q}^*) > \hat{f}(\mathbf{x}, \widehat{\mathbf{Q}}) \Big \}\Big) \leq -r \quad \forall \mathbf{x} \in X, \; \forall \mathbb{Q}^* \in \mathcal{Q} \nonumber
	\end{align}
\end{subequations}
Formally, in the vector optimization problem (\ref{prediction problem}) we aim at finding the least conservative prediction rule~(with respect to the partial order defined above) that satisfies the asymptotic guarantee (\ref{eq: asymptotic guarantee}).

Analogously, we define a partial order on the space of prescription rules and introduce a prescription problem. That is, a prediction-prescription pair $(\hat{f}_1, \hat{x}_1)$ not strictly dominates a pair $(\hat{f}_2, \hat{x}_2)$, if the following set of inequalities holds:
\begin{equation} \label{eq: partial order prescription}
\hat{f}_1(\hat{x}_1(\mathbf{Q}), \mathbf{Q}) \leq \hat{f}_2(\hat{x}_2(\mathbf{Q}), \mathbf{Q}) \quad \forall \mathbb{Q} \in \mathcal{Q},
\end{equation}
\looseness-1where we recall that $\mathbf{Q} = (\mathbb{Q}_1, \ldots, \mathbb{Q}_n)^\top$ is a vector of marginal distributions induced by some $\mathbb{Q} \in \mathcal{Q}$. Let $\mathcal{X}$ be a set of all pairs $(\hat{f}, \hat{x})$, where $\hat{f} \in \mathcal{F}$ and $\hat{x}$ is the prescription rule induced by $\hat{f}$. Then a \textit{prescription problem} can be formulated as follows:
\allowdisplaybreaks\begin{subequations}
	\begin{align} \label{prescription problem} \tag{$\mathcal{P}'$}
	& \min_{(\hat{f}, \hat{x}) \in \mathcal{X}} (\hat{f}, \hat{x}) \quad (\mbox{ w.r.t. }\preceq) \\
	\mbox{\upshape s.t.}
	& \lim_{T_{min} \rightarrow +\infty} \sup \frac{1}{T_{min}}\ln\Big( \Pr\Big\{f(\hat{x}(\widehat{\mathbf{Q}}), \mathbb{Q}^*) > \hat{f}(\hat{x}(\widehat{\mathbf{Q}}), \widehat{\mathbf{Q}}) \Big \}\Big) \leq -r \quad \forall \mathbb{Q}^* \in \mathcal{Q} \nonumber
	\end{align}
\end{subequations}

\looseness-1Table \ref{tab: difference} outlines some comparison of the prediction and prescription problems, (\ref{prediction problem}) and (\ref{prescription problem}), with the related meta-optimization problems in \cite{VanParys2020}. In fact, the key difference is stipulated by our construction of the data set (\ref{data set}) and the objective criterion in (\ref{stochastic programming problem}); recall Assumption \textbf{A3}. Furthermore, in view of Definition \ref{def: data-driven predictor and prescriptor}, our prediction and prescription rules are modeled as some functions of marginal distributions, but not a joint distribution of the cost vector $\mathbf{c}$. 
In contrast to \cite{VanParys2020}, we cannot guarantee that the vector optimization problems (\ref{prediction problem}) and~(\ref{prescription problem}) admit a unique strongly optimal solution. However in the next section we demonstrate that both problems (\ref{prediction problem}) and (\ref{prescription problem}) admit a particular class of weakly optimal solutions that possess some properties of a strongly optimal solution. 

\begin{table}
	\footnotesize
	\centering
	\onehalfspacing
	\begin{tabular}{c | c | c}
		\hline
		Properties & Our formulation & Formulation of Van~Parys~et~al.~\cite{VanParys2020}\\
		\hline
		\multirow{2}{*}{Support set} & \multirow{2}{*}{a discrete set for each component $c_a$, $a \in \mathcal{A}$} & a discrete or continuous set \\
		& & for the cost vector $\mathbf{c}$ \\\hline
		\multirow{2}{*}{Sample size} & \multirow{2}{*}{$T_a \in \mathbb{Z}_{>0}$ for each $a \in \mathcal{A}$} & \multirow{2}{*}{the same size $T \in \mathbb{Z}_{>0}$ for each $a \in \mathcal{A}$}\\
		& & \\\hline
		\multirow{2}{*}{Set of feasible decisions} & \multirow{2}{*}{ $X$ is a bounded subset of $\mathbb{Z}^{n_1}_+ \times \mathbb{R}^{n_2}_+$} & \multirow{2}{*}{$X$ is bounded subset of $\mathbb{R}^{n}$} \\
		& & \\\hline
		\multirow{2}{*}{Objective function in (\ref{stochastic programming problem})} & \multirow{2}{*}{$\gamma(\mathbf{c}, \mathbf{x}) = \mathbf{c}^\top \mathbf{x}$} & $\gamma(\mathbf{c}, \mathbf{x})$ is any continuous function \\
		& & with respect to $\mathbf{x} \in X$ \\\hline
		\multirow{3}{*}{Prediction rules} & an arbitrary function & \multirow{2}{*}{a continuous function} \\
		& of decisions and & \multirow{2}{*}{of decisions and a joint distribution} \\
		&  marginal distributions 
		& \\
		\hline
	\end{tabular}
	\caption{\scriptsize The table summarizes the difference between our formulations (\ref{prediction problem}) and (\ref{prescription problem}) and the related formulations in \cite{VanParys2020}.}
	\label{tab: difference}
\end{table}

 Despite the fact that our approach yields some simplifications compared to the approach of Van~Parys~et~al.~\cite{VanParys2020}, we may capture a number of discrete optimization problem settings, for which an application of the model in \cite{VanParys2020} is substantially limited; recall our discussion of the research questions \textbf{Q1} and \textbf{Q2} in Section \ref{subsec: approach and contribution}. As a byproduct, our construction of the training data set gives rise to some practical insights concerning a relation between the sample size for each component of the cost vector and the form of the nominal distribution; see the research question \textbf{Q3} and Section \ref{sec: comp study} for further~details.



\section{Solution techniques} \label{sec: solution techniques}
In this section we analyze a particular class of distributionally robust prediction and prescription rules based on a relative entropy distance from empirical marginal distributions. More precisely, for any fixed $a \in \mathcal{A}$ the relative entropy distance or Kullback--Leibler divergence \cite{Kullback1997} between an empirical distribution $\widehat{\mathbb{Q}}_{a, T_a}$ and some other distribution $\mathbb{Q}_a \in \mathcal{Q}_a$ is defined as follows:


\begin{equation} \nonumber
D_{KL}(\widehat{\mathbb{Q}}_{a, T_a}\; \Vert \; \mathbb{Q}_a) = \sum_{i = 1}^{d_a} \widehat{q}_{a,i} \ln\frac{\widehat{q}_{a,i}}{q_{a,i}}
\end{equation}

 The idea is to account all marginal distributions within relative entropy balls of some radius $r_a \in \mathbb{R}_{>0}$ centered at the empirical marginal distributions $\widehat{\mathbb{Q}}_{a, T_a}$, $a \in \mathcal{A}$. Specifically, we consider the following component-wise and joint ambiguity sets:
\begin{align} 
& \widehat{\mathcal{Q}}_a := \Big\{ \mathbb{Q}_a \in \mathcal{Q}_0(\mathcal{S}_a): \; 0 \leq q_{a,i} \leq 1 \quad \forall i \in \{1, \ldots, d_a\}, \; \sum_{i = 1}^{d_a} q_{a,i} = 1, \; D_{KL}(\widehat{\mathbb{Q}}_{a,T_a}\; \Vert \; \mathbb{Q}_a) \leq r_a \Big\} \label{empirical marginal ambiguity set} \\
& \widehat{\mathcal{Q}} := \Big\{ \mathbb{Q} \in \mathcal{Q}_0(\mathcal{S}_1 \times \ldots \times \mathcal{S}_n): \mathbb{Q}_a \in \widehat{\mathcal{Q}}_a \quad \forall a \in \mathcal{A} \Big\} \label{ambiguity set}
\end{align}
\looseness-1 The parameters $r_a$, $a \in \mathcal{A}$, are used to control the exponential decay rate induced by the asymptotic guarantee (\ref{eq: asymptotic guarantee}); see Theorem \ref{theorem 1} and Section \ref{subsec: improved finite sample} for further details. 

Based on the definitions above, for some fixed $\mathbf{x} \in X$ we may define the following distributionally robust prediction and prescription rules:
\begin{subequations}
\begin{align}
& \hat{f}_{DR}(\mathbf{x}, \widehat{\mathbf{Q}}) = \max_{\mathbb{Q} \in \widehat{\mathcal{Q}}}\mathbb{E}_{\mathbb{Q}}\{\mathbf{c}^\top \mathbf{x}\} = \sum_{a \in \mathcal{A}} \Big(\max_{\mathbb{Q}_a \in \widehat{\mathcal{Q}}_a}\mathbb{E}_{\mathbb{Q}_a}\{c_a\}\Big) x_a \label{distributionally robust predictor} \\
& \hat{x}_{DR}(\widehat{\mathbf{Q}}) \in \argmin_{\mathbf{x} \in X} \hat{f}_{DR}(\mathbf{x}, \widehat{\mathbf{Q}})
\label{distributionally robust prescriptor}
\end{align}
\end{subequations}

\looseness-1 The remainder of this section is organized as follows. In Sections \ref{subsec: predictor} and \ref{subsec: prescriptor} we establish some theoretical properties of (\ref{distributionally robust predictor}) and (\ref{distributionally robust prescriptor}) in relation to the prediction and prescription problems (\ref{prediction problem}) and (\ref{prescription problem}), respectively. In particular, we introduce a dual reformulation of the maximization problem in (\ref{distributionally robust predictor}) and thereby obtain a linear MIP reformulation of the DRO problem in (\ref{distributionally robust prescriptor}). Finally, Section \ref{subsec: improved finite sample} provides some additional discussion on the choice of the parameters $r_a$, $a \in \mathcal{A}$, in the ambiguity set~(\ref{ambiguity set}). 

\subsection{Properties of the distributionally robust prediction rule} \label{subsec: predictor}
\subsubsection{Dual reformulation} \label{subsec: dual reformulation}
In this section we demonstrate that optimization problem (\ref{distributionally robust predictor}) admits a dual reformulation. In addition, the related DRO problem
\begin{equation} \label{DRO problem} \tag{$\mathcal{DRO}$}
\min_{\mathbf{x} \in X} \hat{f}_{DR}(\mathbf{x}, \widehat{\mathbf{Q}}) = \min_{\mathbf{x} \in X} \max_{\mathbb{Q} \in \widehat{\mathcal{Q}}} \mathbb{E}_{\mathbb{Q}}\{\mathbf{c}^\top \mathbf{x}\}
\end{equation}
can be solved by leveraging $n$ one-dimensional convex optimization problems and a single deterministic MIP problem. The following result holds. 
\begin{proposition} \label{proposition 1}
 The distributionally robust optimization problem (\ref{DRO problem}) is equivalent to the following minimization problem:
\begin{equation} \label{MIP reformulation}
\min_{\mathbf{x} \in X} \sum_{a \in \mathcal{A}} \min_{\beta_a \geq \overline{z}_a} \Big( \beta_a - e^{-r_a}\prod_{i = 1}^{d_a} (\beta_a - z_{a,i})^{\widehat{q}_{a,i}} \Big) x_a,
\end{equation}
where $\overline{z}_a = \max_{i \in \{1, \ldots d_a\}} z_{a,i}$ for each $a \in \mathcal{A}$.	
\begin{proof}
Note that a min-min reformulation of (\ref{DRO problem}) can be obtained by dualizing the worst-case expectation problems for each particular $a \in \mathcal{A}$, that is,
\begin{subequations} \label{worst-case expectation problem}
\begin{align}
& \max_{\mathbf{q}_a}\mathbb{E}_{\mathbb{Q}_a}\{c_a\} = \sum_{i = 1}^{d_a} z_{a,i}q_{a,i} \\
\mbox{s.t. } & 0 \leq q_{a,i} \leq 1 \quad \forall i \in \{1, \ldots, d_a\} \\
& \sum_{i = 1}^{d_a} q_{a,i} = 1 \\
& \sum_{i = 1}^{d_a} \widehat{q}_{a,i} \ln\frac{\widehat{q}_{a,i}}{q_{a,i}} \leq r_a
\end{align}
\end{subequations}	
It is rather straightforward to verify that the feasible set in (\ref{worst-case expectation problem}) is convex and the objective function is linear. By Proposition 2 in \cite{VanParys2020} the dual problem can be reduced to a one-dimensional convex optimization problem of the form:
\begin{subequations} \label{individual dual problems}
\begin{align}
& \min_{\beta_a} \Big( \beta_a - e^{-r_a}\prod_{i = 1}^{d_a} (\beta_a - z_{a,i})^{\widehat{q}_{a,i}} \Big) \\
\mbox{s.t. } & \beta_a \geq \overline{z}_a
\end{align}
\end{subequations}	
If we combine the minimization over $\mathbf{x} \in X$ and the minimization over the dual variables $\beta_a$, $a \in \mathcal{A}$, then the result follows.  		
\end{proof}
\end{proposition}

 Proposition \ref{proposition 1} provides a two-stage solution approach for the DRO problem (\ref{DRO problem}); recall the research question \textbf{Q2}. In the first stage we obtain the worst-case expected costs for each $a \in \mathcal{A}$ by solving univariate convex optimization problems (\ref{individual dual problems}), and in the second stage we solve an instance of the deterministic MIP problem induced by (\ref{MIP reformulation}). Despite the fact that the feasible set in (\ref{individual dual problems}) is unbounded, it can be argued that the optimal solution is unique; we refer to the proof of Proposition~2 in \cite{VanParys2020} for more details. 

\subsubsection{Asymptotic and finite sample guarantees} \label{subsec: asymptotic and finite sample guarantees}
\looseness-1 It turns out that the distributionally robust prediction rule (\ref{distributionally robust predictor}) with an appropriate choice of the parameters $r_a$, $a \in \mathcal{A}$, satisfies both finite sample and asymptotic guarantees, (\ref{eq: finite sample guarantee}) and (\ref{eq: asymptotic guarantee}). In particular, this assertion follows from a strong \textit{large deviation principle}; see, e.g., \cite{Groeneboom1979, VanParys2020}, which is formulated and briefly discussed below.

\begin{proposition} [\textbf{Strong LDP, univariate case}] \label{proposition 2} 
Assume that $a \in \mathcal{A}$ is fixed and the samples $\widehat{c}_{a,1}, \ldots, \widehat{c}_{a,T_a}$ are drawn independently from some marginal distribution $\mathbb{Q}^*_a \in \mathcal{Q}_a$, where $\mathcal{Q}_a$ is given by equation (\ref{marginal ambiguity set}). Then for every Borel set $\mathcal{D}_a \subseteq \mathcal{Q}_a$ and any $T_a \in \mathbb{N}$ we have:	
	
\begin{equation} \label{eq: strong LDP 1} \tag{$\mathcal{LDP}_1$}
\mbox{\upshape Pr} \Big\{\widehat{\mathbb{Q}}_{a,T_a} \in \mathcal{D}_a \Big\} \leq (T_a + 1)^{d_a}\exp \Big(-T_a \inf_{\mathbb{Q}_a \in \mathcal{D}_a} D_{KL}(\mathbb{Q}_a \; \Vert \; \mathbb{Q}^*_a)\Big)
\end{equation}
\begin{proof} See Theorem 2 in \cite{VanParys2020}. 
\end{proof}
\end{proposition}

In other words, if $\mathbb{Q}^*_a \notin \mathcal{D}_a$, then the empirical distribution $\widehat{\mathbb{Q}}_{a,T_a}$ belongs to $\mathcal{D}_a$ with a probability that decays exponentially in the sample size, $T_a$. In the next result we mimic the idea of Example~\ref{example 1}. That is, the out-of-sample prediction disappointment (\ref{eq: out-of-sample prediction}) with respect to (\ref{distributionally robust predictor}) is bounded from above by some function that decays exponentially in $T_{min}$.

\begin{theorem} \label{theorem 1}
Assume that $r > 0$ is a required exponential decay rate, $T_{min} = \min_{a \in \mathcal{A}}T_a$ and $\delta_a \in (0, 1)$ are such that $\sum_{a \in \mathcal{A}}\delta_a = 1$. Then the prediction rule (\ref{distributionally robust predictor}) with the parameters 
\begin{equation} \label{eq: radius}
r_a = \frac{1}{T_a}\Big(d_a \ln (T_a + 1) + r T_{min} - \ln \delta_a\Big) \quad \forall a \in \mathcal{A}
\end{equation}
is feasible in (\ref{prediction problem}). 

\begin{proof}
	
We observe that $\hat{f}_{DR} \in \mathcal{F}$ by construction and, thus, it suffices to prove that $\hat{f}_{DR}$ satisfies the asymptotic guarantee (\ref{eq: asymptotic guarantee}). 		
Let $\mathcal{A}_{\mathbf{x}} = \{a \in \mathcal{A}: x_a > 0\} \neq \varnothing$ for each $\mathbf{x} \in X$. We observe that the out-of-sample prediction disappointment (\ref{eq: out-of-sample prediction}) is bounded from above as:
\begin{equation} \label{eq: feasibility proof 1} \nonumber
\begin{gathered}
\Pr\Big\{f(\mathbf{x}, \mathbb{Q}^*) > \hat{f}_{DR}(\mathbf{x}, \widehat{\mathbf{Q}}) \Big \} = \Pr\Big\{\mathbb{E}_{\mathbb{Q}^*}\{\mathbf{c}^\top \mathbf{x}\} - \max_{\mathbb{Q} \in \widehat{\mathcal{Q}}}\mathbb{E}_{\mathbb{Q}}\{\mathbf{c}^\top \mathbf{x}\} > 0\Big\} = \\ \Pr\Big\{\sum_{a \in \mathcal{A}} \Big(\mathbb{E}_{\mathbb{Q}^*_a}\{c_a\} - \max_{\mathbb{Q}_a \in \widehat{\mathcal{Q}}_a}\mathbb{E}_{\mathbb{Q}_a}\{c_a\}\Big)x_a > 0 \Big\} \leq \Pr\Big\{\bigvee_{a \in \mathcal{A}_{\mathbf{x}}} \Big(\mathbb{E}_{\mathbb{Q}^*_a}\{c_a\} - \max_{\mathbb{Q}_a \in \widehat{\mathcal{Q}}_a}\mathbb{E}_{\mathbb{Q}_a}\{c_a\} > 0 \Big) \Big\},
\end{gathered}
\end{equation}
where the last inequality exploits Assumption \textbf{A2} and can be checked by contradiction. Furthermore, the inequality
$$\mathbb{E}_{\mathbb{Q}^*_a}\{c_a\} - \max_{\mathbb{Q}_a \in \widehat{\mathcal{Q}}_a}\mathbb{E}_{\mathbb{Q}_a}\{c_a\} > 0$$
for some $a \in \mathcal{A}$ implies that $\mathbb{Q}^*_a \notin \widehat{\mathcal{Q}}_a$. Hence, by leveraging the union bound and Proposition \ref{proposition 2} we observe that
\begin{equation} \label{eq: feasibility proof 2} \nonumber
\begin{gathered}
\Pr\Big\{\bigvee_{a \in \mathcal{A}_{\mathbf{x}}} \Big(\mathbb{E}_{\mathbb{Q}^*_a}\{c_a\} - \max_{\mathbb{Q}_a \in \widehat{\mathcal{Q}}_a}\mathbb{E}_{\mathbb{Q}_a}\{c_a\} > 0 \Big) \Big\} \leq \sum_{a \in \mathcal{A}_{\mathbf{x}}} \Pr\Big\{ \mathbb{E}_{\mathbb{Q}^*_a}\{c_a\} - \max_{\mathbb{Q}_a \in \widehat{\mathcal{Q}}_a}\mathbb{E}_{\mathbb{Q}_a}\{c_a\} > 0 \Big\} \leq \\ \sum_{a \in \mathcal{A}_{\mathbf{x}}} \Pr\Big\{\mathbb{Q}^*_a \notin \widehat{\mathcal{Q}}_a\Big\} = \sum_{a \in \mathcal{A}_{\mathbf{x}}} \Pr\Big\{D_{KL}(\widehat{\mathbb{Q}}_{a,T_a} \; \Vert \; \mathbb{Q}^*_a) > r_a\Big\} \leq \sum_{a \in \mathcal{A}_{\mathbf{x}}}(T_a + 1)^{d_a} e^{-T_a r_a} 
\end{gathered}
\end{equation}
 
 By leveraging equation (\ref{eq: radius}) we conclude that
 \begin{equation} \nonumber
 \sum_{a \in \mathcal{A}_{\mathbf{x}}}(T_a + 1)^{d_a} e^{-T_a r_a} = \sum_{a \in \mathcal{A}_{\mathbf{x}}} \delta_a e^{-T_{min} r} \leq e^{-T_{min}r}
 \end{equation} 
 Hence,
 \begin{equation} \nonumber
 \begin{gathered}
 \limsup_{T_{min} \to \infty} \frac{1}{T_{min}}\ln\Big(\Pr\Big\{f(\mathbf{x}, \mathbb{Q}^*) > \hat{f}_{DR}(\mathbf{x}, \widehat{\mathbf{Q}}) \Big \}\Big) \leq \limsup_{T_{min} \to \infty} \frac{1}{T_{min}}\ln(e^{-T_{min} r}) = -r
 \end{gathered}
 \end{equation}
and the proof follows.
\end{proof}
\end{theorem}

We infer that the distributionally robust prediction rule (\ref{distributionally robust predictor}) with the parameters $r_a$, $a \in \mathcal{A}$, given by equation (\ref{eq: radius}) is feasible in the prediction problem (\ref{prediction problem}). In the next section we demonstrate that (\ref{distributionally robust predictor}) is, in fact, a weakly optimal solution of (\ref{prediction problem}). 

\subsubsection{Weak optimality} \label{subsubsec: weak optimality}
At first, we introduce some problem-specific definitions of weak and strong optimality. In particular, we introduce a class of weakly optimal solutions that dominate all other feasible solutions at some given set of points. 
 
\begin{definition}[\textbf{Weak optimality, prediction}] \label{def: weak optimality} \upshape
A prediction rule $\hat{f} \in \mathcal{F}$ is called \textit{weakly optimal} (\textit{W}-optimal) for (\ref{prediction problem}), if it is feasible in (\ref{prediction problem}) and any other prediction rule $\hat{f}' \in \mathcal{F}$ that satisfies $$\hat{f}'(\mathbf{x}, \mathbf{Q}) < \hat{f}(\mathbf{x}, \mathbf{Q}) \quad \forall \mathbf{x} \in X, \; \forall \mathbb{Q} \in \mathcal{Q}$$ is infeasible in (\ref{prediction problem}). Furthermore, $\hat{f} \in \mathcal{F}$ is \textit{strongly optimal with respect to a nonempty set $\widetilde{\mathcal{Z}} \subseteq X \times \mathcal{Q}$} (\textit{S}-optimal w.r.t. $\widetilde{\mathcal{Z}}$), if it is feasible in (\ref{prediction problem}) and any other \textit{feasible} prediction rule, $\hat{f}' \in \mathcal{F}$, satisfies 
$$\hat{f}'(\mathbf{x}, \mathbf{Q}) \geq \hat{f}(\mathbf{x}, \mathbf{Q}) \quad \forall (\mathbf{x}, \mathbb{Q}) \in \widetilde{\mathcal{Z}} $$
\vspace{-9.5mm}\flushright$\square$ 
\end{definition}

As a remark, if $\widetilde{\mathcal{Z}} = X \times \mathcal{Q}$ in the definition above, then the related prediction rule is strongly optimal in a common sense. 
It can be also argued that an \textit{S}-optimal w.r.t $\widetilde{\mathcal{Z}} \neq \varnothing$ prediction rule is also \textit{W}-optimal, but not vice versa. We illustrate this fact with the following counterexample.

\begin{example} \label{example 2} \upshape
We consider the following vector optimization problem:
\begin{align} \label{toy vector optimization problem}
& \min_{\phi(\cdot)} \phi(y) \\
\mbox{s.t. } & \int_{0}^{1}\phi(y)dy \geq \frac{1}{2} \nonumber \\
& \phi(0) \geq 0 \nonumber
\end{align}	
where $\phi: [0, 1] \rightarrow \mathbb{R}$ is an integrable function on the interval $[0,1]$. It is also assumed that $\phi_1 \preceq \phi_2$ if and only if $\phi_1(y) \leq \phi_2(y)$ for all $y \in [0, 1]$. 

\looseness-1 We note the function $\phi_1(y) := 0.5$ is a \textit{W}-optimal solution of (\ref{toy vector optimization problem}) as the function $\widetilde{\phi}_1(y) := 0.5 - \varepsilon$ does not satisfy the first constraint in (\ref{toy vector optimization problem}) for any $\varepsilon > 0$. However this solution is manifestly not \textit{S}-optimal w.r.t. any nonempty set $\widetilde{\mathcal{Z}} \subseteq [0, 1]$, as the function $\phi_1(y)$ can be decreased at any fixed point $y_0 \in [0,1]$ without changing the value of the integral in (\ref{toy vector optimization problem}). 
On the other hand, the function $\phi_2(y) := y$ is an \textit{S}-optimal w.r.t. $\widetilde{\mathcal{Z}} = \{0\}$ solution of~(\ref{toy vector optimization problem}) as any function $\widetilde{\phi}_2(y)$ such that $\widetilde{\phi}_2(0) < \phi_2(0) = 0$ is infeasible in (\ref{toy vector optimization problem}). 

In summary, we observe that the functions $\phi_1(y)$ and $\phi_2(y)$ are not comparable in general. However the function $\phi_2(y)$ can be posed as the best solution at least at one point, i.e., $y = 0$. One possible interpretation is that the more ``rich'' the set $\widetilde{\mathcal{Z}}$, the more advantage an $S$-optimal w.r.t. $\widetilde{\mathcal{Z}}$ solution has compared to some typical $W$-optimal solutions.  \vspace{-9.5mm} \flushright $\square$
\end{example}

It turns out that \textit{W}-optimality holds for a rather broad class of prediction rules including the prediction rule (\ref{eq: sample-average predictor + epsilon}) from Example \ref{example 1} and the distributionally robust prediction rule (\ref{distributionally robust predictor}). 
The following result takes place. 

\begin{proposition} \label{proposition 3}
	Let $\overline{z}_a = \max_{i \in \{1, \ldots, d_a\}} z_{a,i}$ for any $a \in \mathcal{A}$. Then any prediction rule $\hat{f} \in \mathcal{F}$ that is feasible in (\ref{prediction problem}) and satisfies the following upper bound:
	\begin{equation} \label{eq: WW optimality upper bound}
	\hat{f}(\mathbf{x}, \mathbf{Q}) \leq \sum_{a \in \mathcal{A}} \overline{z}_a x_a \quad \forall \mathbf{x} \in X, \; \forall \mathbb{Q} \in \mathcal{Q}
	\end{equation}
	is W-optimal in (\ref{prediction problem}). 	
	\begin{proof} See \ref{sec: app}. 
	\end{proof}	
\end{proposition}

 Proposition \ref{proposition 3} does not exploit the structure of $\hat{f}$, but only requires that the estimate of expected loss, $\hat{f}(\mathbf{x}, \mathbf{Q})$, satisfies the asymptotic guarantee (\ref{eq: asymptotic guarantee}) and does not exceed the maximal possible loss; see the right hand side of equation (\ref{eq: WW optimality upper bound}). In the following we demonstrate that the distributionally robust prediction rule (\ref{distributionally robust predictor}) is also \textit{S}-optimal w.r.t. to some specified set $\widetilde{\mathcal{Z}} \subseteq X \times \mathcal{Q}$. The subsequent theoretical results exploit the following additional assumption on the individual support sets $\mathcal{S}_a$, $a \in \mathcal{A}$:
\begin{itemize}
	\item[\textbf{A1'.}] For any $a, b \in \mathcal{A}$, $a \neq b$, we have $d_a = d_b$ and there exist some $w_{a, b} \in \mathbb{R}_+$ and $v_{a, b} \in \mathbb{R}$ such that $z_{b, i} = w_{a, b}\; z_{a, i} + v_{a, b}$ for any $i \in \{1, \ldots, d_a\}$.
\end{itemize}

In fact, Assumption \textbf{A1'} entails an affine dependence between the elements of individual support sets. Besides that, the number of elements $d_a$, $a \in \mathcal{A}$, is assumed to be the same for all components of the cost vector $\mathbf{c}$. While Assumption \textbf{A1'} limits the number of candidate data-generating distributions, we may still account different cost ranges for the components of $\mathbf{c}$ by adapting the constants $w_{a, b}$ and $v_{a, b}$ for $a, b \in \mathcal{A}$, $a \neq b$. 
The following result holds.

\begin{theorem} \label{theorem 2}
	Assume that $r > 0$ is a required exponential decay rate and Assumption \textbf{A1'} holds. Also, let
	\begin{equation} \label{eq: set Z}
	\widetilde{\mathcal{Z}} = X \times \Big\{\mathbb{Q} \in \mathcal{Q}: q_{a,i} = q_{b,i} > 0 \quad \forall i \in \{1, \ldots, d_a\}, \; \forall a,b \in \mathcal{A}, \; a \neq b
	 \Big\},
	\end{equation}
  where $q_{a,i}$ are rational and denote the probability that $c_a$ equals $z_{a,i}$. 	
	 Then the distributionally robust prediction rule (\ref{distributionally robust predictor}) with parameters $r_a$, $a \in \mathcal{A}$, given by equation (\ref{eq: radius}) with $\delta_a = \frac{1}{n}$ is \textit{S}-optimal w.r.t. $\widetilde{\mathcal{Z}}$ for (\ref{prediction problem}).
	\begin{proof}

	We need to show for any fixed $(\widetilde{\mathbf{x}}, \widetilde{\mathbb{Q}}) \in \widetilde{\mathcal{Z}}$ every other prediction rule $\hat{f}' \in \mathcal{F}$ that satisfies
	\begin{equation} \label{eq: new weak optimality proof 1}
	\varepsilon := \hat{f}_{DR}(\widetilde{\mathbf{x}}, \widetilde{\mathbf{Q}}) - \hat{f}'(\widetilde{\mathbf{x}}, \widetilde{\mathbf{Q}}) > 0 
	\end{equation}
	is infeasible in (\ref{prediction problem}); recall the definition of \textit{S}-optimality. 	
	
	
	\noindent
	\textbf{Step 1.} At the first step bound the out-of-sample disappointment 
	\begin{equation} \nonumber
	\Pr\Big\{f(\mathbf{x}, \mathbb{Q}^*) > \hat{f}'(\mathbf{x}, \widehat{\mathbf{Q}}) \Big\}
	\end{equation} 
	from below at $\mathbf{x} = \widetilde{\mathbf{x}}$ by selecting some marginal probability distributions $\mathbb{Q}^*_a$ for each $a \in \mathcal{A}$. In this regard, for each $a \in \mathcal{A}$ we construct a relative entropy ball $\widehat{\mathcal{Q}}_a$ of the radius $r_a$ centered at the marginal distribution $\widetilde{\mathbb{Q}}_a$; see equation (\ref{empirical marginal ambiguity set}) and Figure~\ref{fig: relative entropy ball} for details. In particular, by definition of $\widetilde{\mathcal{Z}}$ the marginal distributions $\widetilde{\mathbb{Q}}_a$, $a \in \mathcal{A}$, are the same probability distributions but equipped with different support sets $\mathcal{S}_a$.
	
	Then we set $b \in \argmin_{a \in \mathcal{A}}T_a$ and define the worst-case marginal distributions as:
	\begin{equation} \nonumber
	\mathbb{Q}_a^{(w)} \in \argmax_{\mathbb{Q}_a \in \widehat{\mathcal{Q}}_a}\mathbb{E}_{\mathbb{Q}_a}\{c_a\}
	\end{equation}
	In view of Assumption \textbf{A1'}, we have $d_a = d \in \mathbb{Z}_{>0}$ for each $a \in \mathcal{A}$.
	Furthermore, substituting $\delta_a = \frac{1}{n}$, $a \in \mathcal{A}$, into equation (\ref{eq: radius}) implies that
	$$r_a = \frac{1}{T_a}\Big(d \ln (T_a + 1) + r T_{min} - \ln \frac{1}{n}\Big) \quad \forall a \in \mathcal{A}$$
  Hence, there exists $T_0 \in \mathbb{Z}_{>0}$ such that for every $T_{min} = T_b \geq T_0$ $$r_b \geq r_a \quad \forall a \in \mathcal{A},\; a \neq b$$
	Then using Assumption \textbf{A1'} for any fixed $a \in \mathcal{A}$, $a \neq b$, and $T_{min} \geq T_0$ we observe that
	\begin{equation} \label{eq: weak optimality proof A1'}
	\begin{gathered}
	 \mathbb{E}_{\mathbb{Q}^{(w)}_a} \{c_a\} = \sum_{i = 1}^{d} q^{(w)}_{a, i}z_{a,i} = \max_{\mathbb{Q}_a \in \widehat{\mathcal{Q}}_a} \sum_{i = 1}^{d} q_{a, i}z_{a,i} = \max_{\mathbb{Q}_a \in \widehat{\mathcal{Q}}_a} \sum_{i = 1}^{d} q_{a, i}(w_{b, a} z_{b, i} + v_{b, a}) = \\ w_{b, a} \max_{\mathbb{Q}_a \in \widehat{\mathcal{Q}}_a} \Big( \sum_{i = 1}^{d} q_{a, i} z_{b, i} \Big) + v_{b, a} \leq w_{b, a} \max_{\mathbb{Q}_b \in \widehat{\mathcal{Q}}_b} \Big( \sum_{i = 1}^{d} q_{b, i} z_{b, i} \Big) + v_{b, a} = \\ \max_{\mathbb{Q}_b \in \widehat{\mathcal{Q}}_b} \Big( \sum_{i = 1}^{d} q_{b, i} (w_{b, a} z_{b, i} + v_{b, a})\Big) = \mathbb{E}_{\mathbb{Q}^{(w)}_b} \{c_a\} 
	\end{gathered}
	\end{equation} 
	The inequality follows from the fact that $w_{b,a} > 0$ and the relative entropy balls $\widehat{\mathcal{Q}}_a$ and $\widehat{\mathcal{Q}}_b$ have the same center (recall the definition of $\widetilde{\mathcal{Z}}$).
	

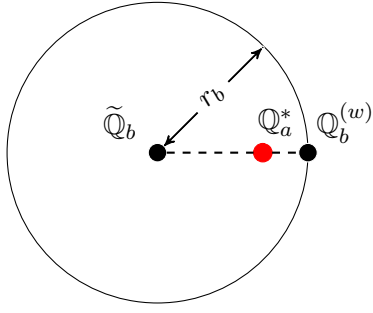
\begin{figure}
	\centering
	\begin{tikzpicture}[scale=1,transform shape]
	\filldraw[fill=white, draw=black] (0,0) circle (2cm);
	\tikzset{VertexStyle/.style = {draw=white,fill=black,circle,
			minimum size=7pt,inner sep=0pt}}  
	\Vertex[x=0,y=0,L=$ $]{1}
	\Vertex[x=2,y=0,L=$ $]{4}	
	\tikzset{VertexStyle/.style = {draw=white,fill=black,circle,
			minimum size=7pt, color=red, inner sep=0pt}}
		\Vertex[x=1.4,y=0,L=$ $]{3}
	\tikzset{VertexStyle/.style = {draw=white,fill=black,circle,
			minimum size=1pt,inner sep=0pt}} 
	\Vertex[x=1.4,y=1.4,L=$ $]{2}		
	\tikzset{VertexStyle/.style = {draw=white,circle,
			minimum size=1pt,inner sep=0pt}}
	\Vertex[x=-0.5,y=0.4,L=$\widetilde{\mathbb{Q}}_b$]{10}	\Vertex[x=2.5,y=0.4,L=$\mathbb{Q}_b^{(w)}$]{40}
	\Vertex[x=1.57,y=0.4,L=$\mathbb{Q}^*_a$]{50}
	\tikzstyle{LabelStyle}=[fill=white,sloped]
	\tikzstyle{EdgeStyle}=[post]
	\Edge(2)(1)
	\Edge[label=$r_b$](1)(2)
	\tikzstyle{EdgeStyle}=[dashed]
	\Edge(3)(4)
	\Edge(1)(3)
	\end{tikzpicture}
	\caption{\footnotesize A rough illustration of the notations used in the first step of the proof of Theorem \ref{theorem 3}. The relative entropy ball $\widehat{\mathcal{Q}}_b$ is depicted as an Euclidean ball of the radius $r_b$. The red dot corresponds to the nominal marginal distribution $\mathbb{Q}^*_a$, which is the same for all $a \in \mathcal{A}$.}
	\label{fig: relative entropy ball}
\end{figure}

	Next, we define 
	\begin{equation} \label{eq: weak optimality proof nominal marginal distributuins}
	\mathbb{Q}^*_a := \lambda \widetilde{\mathbb{Q}}_b + (1 - \lambda)\mathbb{Q}^{(w)}_b \quad \forall a \in \mathcal{A}
	\end{equation} 
	for some $\lambda \in (0, 1)$. By definition of $\widetilde{\mathcal{Z}}$ we have $\widetilde{q}_{b,i} > 0$, $i \in \{1, \ldots, d\}$, and, thus, the nominal probabilities $q^*_{a,i}$ are also positive and independent of $a \in \mathcal{A}$. Consequently, the relative entropy distance $D_{KL}(\widetilde{\mathbb{Q}}_b \; \Vert \; \mathbb{Q}^*_a)$ is finite. 
	
	We observe that convexity of the relative entropy distance implies that:
	\begin{equation} \label{eq: new weak optimality proof 2.1}
	D_{KL}(\widetilde{\mathbb{Q}}_b \; \Vert \; \mathbb{Q}^*_a) \leq \lambda D_{KL}(\widetilde{\mathbb{Q}}_b \; \Vert \; \widetilde{\mathbb{Q}}_b) + (1 - \lambda) D_{KL}(\widetilde{\mathbb{Q}}_b \; \Vert \; \mathbb{Q}^{(w)}_b) \leq (1 - \lambda)r_b < r_b
	\end{equation}
	Furthermore, for any $\varepsilon' > 0$ and $a \in \mathcal{A}$ we may guarantee that for $T_{min} \geq T_0$
	\begin{equation} \label{eq: new weak optimality proof 2.2}
	\mathbb{E}_{\mathbb{Q}^*_a}\{c_a\} = \lambda \mathbb{E}_{\widetilde{\mathbb{Q}}_b}\{c_a\} + (1 - \lambda)\mathbb{E}_{\mathbb{Q}^{(w)}_b}\{c_a\} > \mathbb{E}_{\mathbb{Q}^{(w)}_b}\{c_a\} - \varepsilon' \geq \mathbb{E}_{\mathbb{Q}^{(w)}_a}\{c_a\} - \varepsilon',
	\end{equation} 
	where the strict inequality can be always satisfied with a sufficiently small $\lambda > 0$ and the last inequality exploits (\ref{eq: weak optimality proof A1'}). 
		
	Let $\varepsilon' = \frac{\varepsilon}{nu}$, where $u = \max\{u_a, a \in \mathcal{A}\}$ is a maximal upper bound for the decision variables induced by Assumption \textbf{A2}. Then by leveraging (\ref{eq: new weak optimality proof 1}) and (\ref{eq: new weak optimality proof 2.2}) it can be observed that for $T_{min} \geq T_0$
	\begin{equation} \label{eq: new weak optimality proof 3}
	\begin{gathered}
	f(\widetilde{\mathbf{x}}, \mathbb{Q}^*) = \sum_{a \in \mathcal{A}} \mathbb{E}_{\mathbb{Q}^*_a}\{c_a\}\widetilde{x}_a > \sum_{a \in \mathcal{A}} \Big(\max_{\mathbb{Q}_a \in \widehat{\mathcal{Q}}_a}\mathbb{E}_{\mathbb{Q}_a}\{c_a\} - \varepsilon' \Big)\widetilde{x}_a = \\
	= \hat{f}_{DR}(\widetilde{\mathbf{x}}, \widetilde{\mathbf{Q}}) - \varepsilon'\sum_{a \in \mathcal{A}} \widetilde{x}_a = \hat{f}^{'}(\widetilde{\mathbf{x}}, \widetilde{\mathbf{Q}}) + \varepsilon - \frac{\varepsilon}{nu} \sum_{a \in \mathcal{A}} \widetilde{x}_a \geq \hat{f}^{'}(\widetilde{\mathbf{x}}, \widetilde{\mathbf{Q}})
	\end{gathered} 
	\end{equation} 
	
	In conclusion, we need to show that $\widetilde{\mathbf{Q}}$ can be posed as some vector of \textit{empirical marginal distributions}. In this particular case using the inequality (\ref{eq: new weak optimality proof 3}) we may bound the out-of-sample disappointment from below as: 
	\begin{equation} \label{eq: new weak optimality proof 5} \tag{$LB_1$}
	\Pr\Big\{f(\widetilde{\mathbf{x}}, \mathbb{Q}^*) > \hat{f}'(\widetilde{\mathbf{x}}, \widehat{\mathbf{Q}}) \Big\} \geq \Pr\Big\{\widehat{\mathbf{Q}} = \widetilde{\mathbf{Q}} \Big\}
	\end{equation}	
	
	Indeed, by definition of $\widetilde{\mathcal{Z}}$ the probabilities $\widetilde{q}_{a,i}$ for $i \in \{1, \ldots, d\}$ and any $a \in \mathcal{A}$ are rational. Therefore, the vector $\widetilde{\mathbf{Q}}$ can be posed as vector of the same empirical marginal distributions on $K \in \mathbb{Z}_{>0}$ samples, where $K$ is the smallest common multiple for the denominators of $\widetilde{q}_{a,i}$, $i~\in~\{1, \ldots, d\}$. Furthermore, it is rather straightforward to verify that the same empirical distributions can be obtained from $2K, 3K, \ldots$ samples as well and without loss of generality we may select~$K \geq T_0$. Eventually, since the asymptotic guarantee (\ref{eq: asymptotic guarantee}) must hold for any $T_1, \ldots, T_n \in \mathbb{Z}_{>0}$ as $T_{min} = \min_{a \in \mathcal{A}} T_a$ goes to infinity, we consider the case, in which
	\begin{equation} \label{eq: new weak optimality proof 4}
	T_a \in \mathcal{K} := \{K, 2K, \ldots\}, \quad \forall a \in \mathcal{A}
	\end{equation}
	That is, for any $T_1, \ldots, T_n$ that satisfy equation (\ref{eq: new weak optimality proof 4}) the lower bound (\ref{eq: new weak optimality proof 5}) holds. 
		
	\noindent \textbf{Step 2.} The goal of this step is to construct a joint distribution $\mathbb{Q}^*$ with the marginal distributions defined by equation (\ref{eq: weak optimality proof nominal marginal distributuins}) such that the asymptotic guarantee (\ref{eq: asymptotic guarantee}) for $\hat{f}'$ does not hold at $\mathbf{x} = \widetilde{\mathbf{x}}$. First,
	we assume that the vector of marginal distributions $\widetilde{\mathbf{Q}}$ is induced by some joint distribution $\widetilde{\mathbb{J}} \in \mathcal{Q}$, which is itself an empirical distribution on $T \in \mathcal{K}$ samples. Formally, as the marginal distributions $\widetilde{\mathbb{Q}}_a$, $a \in \mathcal{A}$, coincide, for any set of indexes $(i_1, \ldots, i_n)^\top \in \{1, \ldots, d\}^n$ and any fixed $T \in \mathcal{K}$ we may set
	\begin{equation} \nonumber
	{\Pr}_{\;\widetilde{\mathbb{J}}}\{\mathbf{c} = (z_{1,i_1}, \ldots, z_{n,i_n})^\top\} = j_{i_1, \ldots, i_n} = \begin{cases}
	\widetilde{q}_{1,i}, \mbox { if } i_1 = \ldots = i_n = i \mbox{ for } i \in \{1, \ldots, d\}\\
	0, \mbox{ otherwise}
	\end{cases}
	\end{equation}
	Analogously, since the marginal distributions $\mathbb{Q}^*_a$ coincide for each $a \in \mathcal{A}$, we may define the nominal distribution, $\mathbb{Q}^*$, as:
	\begin{equation} \nonumber
	{\Pr}_{\mathbb{Q}^*}\{\mathbf{c} = (z_{1,i_1}, \ldots, z_{n,i_n})^\top\} = {q}^*_{i_1, \ldots, i_n} = \begin{cases}
	{q}^*_{1,i}, \mbox { if } i_1 = \ldots = i_n = i \mbox{ for } i \in \{1, \ldots, d\} \\
	0, \mbox{ otherwise}
	\end{cases}
	\end{equation}
	
	We conclude that for $T \in \mathcal{K}$ the following lower bound holds:
	\begin{equation} \label{eq: new weak optimality proof 6} \tag{$LB_2$}
	\begin{gathered}
	\Pr\Big\{\widehat{\mathbf{Q}} = \widetilde{\mathbf{Q}} \Big\} \geq \Pr\Big\{\widehat{\mathbb{J}} = \widetilde{\mathbb{J}} \Big\} \geq 
	(T + 1)^{-d^n} e^{-TD_{KL}(\widetilde{\mathbb{J}} \; \Vert \; \mathbb{Q}^*)} = \\
	(T + 1)^{-d^n} e^{-TD_{KL}(\widetilde{\mathbb{Q}}_b \; \Vert \; \mathbb{Q}^*_b)} \geq (T + 1)^{-d^n} e^{-(1 - \lambda)r_b T}
	\end{gathered}
	\end{equation}
	In particular, the first inequality in (\ref{eq: new weak optimality proof 6}) follows from the fact the joint distribution $\widetilde{\mathbb{J}}$ enforces the vector of empirical marginal distributions $\widetilde{\mathbf{Q}}$. The second inequality exploits that:
	\begin{itemize}
		\item[(\textit{i})] the number of sequences of size $T$ that give rise to the same empirical joint distribution is bounded from below by $(T + 1)^{-|\mathcal{S}|}$, where $|\mathcal{S}| = d^n$ is a cardinality of the support of $\mathbf{c}$;
		\item[(\textit{ii})] the probability of each sequence coincides with $e^{-TD_{KL}(\widetilde{\mathbb{J}} \; \Vert \; \mathbb{Q}^*)}$.
	\end{itemize}
 We omit the proof of (\textit{i}) and (\textit{ii}) for brevity and refer to the proof of Theorem 1 in \cite{VanParys2020} for details. The equality in (\ref{eq: new weak optimality proof 6}) is implied by the definition of the relative entropy distance. Finally, the last inequality follows from (\ref{eq: new weak optimality proof 2.1}).
	
By combining the lower bounds (\ref{eq: new weak optimality proof 5}) and (\ref{eq: new weak optimality proof 6}) for $T_{min} := T(k) = kK \to +\infty$, $k \in \{1, 2, \ldots\}$, we observe that
	\begin{equation} \nonumber
	\begin{gathered}
	\limsup_{T_{min} \to +\infty} \frac{1}{T_{min}} \ln \Big(\Pr\Big\{f(\widetilde{\mathbf{x}}, \mathbb{Q}^*) > \hat{f}'(\widetilde{\mathbf{x}}, \widehat{\mathbf{Q}}) \Big\}\Big) \geq \lim_{k \to +\infty} \frac{1}{T(k)} \ln\Big((T(k) + 1)^{-d^n} e^{-(1 - \lambda)r_b T(k)}\Big) = \\
	\lim_{k \to +\infty}\frac{1}{T(k)}\Big(-d^n \ln(T(k) + 1) -(1 - \lambda)\big(d \ln(T(k) + 1) + rT(k) - \ln \frac{1}{n}\big) \Big) = -(1 - \lambda)r > -r
	\end{gathered}
	\end{equation} 
	Thus, the prediction rule $\hat{f}'$ is infeasible in (\ref{prediction problem}) and the result follows.	
	\end{proof}
\end{theorem}

As outlined earlier, the set $\widetilde{\mathcal{Z}}$ defined by equation (\ref{eq: set Z}) contains a class of points $(\mathbf{x}, \mathbb{Q})$, where $\mathbf{x} \in X$ is any feasible decision and $\mathbb{Q} \in \mathcal{Q}$ is a joint distributions with the same marginals. At the same time, the proof of Theorem \ref{theorem 2} cannot be applied directly to the case when the parameters $d_a$, $a \in \mathcal{A}$, are different. 
In this regard, we formulate an auxiliary result, which exploits some relaxed version of Assumption~\textbf{A1'} and thereby holds for any $d_a \in \mathbb{Z}_{>0}$. More specifically, we make the following assumption:

\begin{itemize}
	 \item[\textbf{A1''}] Suppose that the values $z_{a,1}, \ldots, z_{a,d_a}$ are sorted in a decreasing order for each $a \in \mathcal{A}$. Then for any $a, b \in \mathcal{A}$, $a \neq b$, there exist some $w_{a, b} \in \mathbb{R}_+$ and $v_{a, b} \in \mathbb{R}$ such that 
$$z_{b, i} = w_{a, b}\; z_{a, i} + v_{a, b}$$ for any $i \leq d_{min}$, where $d_{min} = \min_{a \in \mathcal{A}} d_a$.
\end{itemize}

The following result holds. 
\begin{theorem} \label{theorem 3}
Assume that $r > 0$ is a required exponential decay rate and Assumption \textbf{A1''} holds. Also, let
\begin{equation} \label{eq: set Z 2}
\widetilde{\mathcal{Z}} = X \times \Big\{\mathbb{Q} \in \mathcal{Q}: q_{a,i} = \begin{cases} q_{b,i} \in \mathbb{R}_{>0} \mbox{, if } i \leq d_{min} \\
0, \mbox{ otherwise} \end{cases} \quad \forall a,b \in \mathcal{A}, \; a \neq b
\Big\},
\end{equation}
where $q_{a,i}$ are rational and denote the probability that $c_a$ equals $z_{a,i}$. 	
 Then the distributionally robust prediction rule (\ref{distributionally robust predictor}) with parameters $r_a$, $a \in \mathcal{A}$, given by equation (\ref{eq: radius}) with $\delta_a = \frac{1}{n}$ is \textit{S}-optimal w.r.t. $\widetilde{\mathcal{Z}}$ for (\ref{prediction problem}).
	\begin{proof} See \ref{sec: app}.
	\end{proof}
\end{theorem}

In conclusion, it is still an open question whether (\ref{distributionally robust predictor}) is a strongly optimal solution for (\ref{prediction problem}) or not. We did not manage to answer this question by using the methodology of Theorem~\ref{theorem 2} and, thus, we leave it as a possible direction of future research.


\subsection{Properties of the distributionally robust prescriptor} \label{subsec: prescriptor}

\looseness-1 In this section we consider analogous properties of the prescription rule (\ref{distributionally robust prescriptor}), which is defined by a pair of functions $(\hat{f}_{DR}, \hat{x}_{DR})$. As a corollary of Theorem \ref{theorem 1}, we observe that the pair $(\hat{f}_{DR}, \hat{x}_{DR})$ is feasible in (\ref{prescription problem}). That is, the asymptotic guarantee (\ref{eq: asymptotic guarantee}) holds for any feasible decision $\mathbf{x} \in X$ and, in particular, for the optimal estimate 
$$\hat{x}_{DR}(\widehat{\mathbf{Q}}) \subseteq \argmin_{\mathbf{x} \in X} \hat{f}_{DR}(\mathbf{x}, \widehat{\mathbf{Q}})$$ 
Next, we slightly modify the definitions of weak and strong optimality in the context of the prescription problem (\ref{prescription problem}). 

\begin{definition}[\textbf{Weak and strong optimality, prescription}] \label{def: weak optimality 2} \upshape
A prediction-prescription pair $(\hat{f}, \hat{x})$ is called \textit{weakly optimal} (\textit{W}-optimal) for (\ref{prescription problem}), if it is feasible in (\ref{prescription problem}) and any other prediction-prescription pair $(\hat{f}', \hat{x}')$ that satisfies $$\hat{f}'(\hat{x}'(\mathbf{Q}), \mathbf{Q}) < \hat{f}(\hat{x}(\mathbf{Q}), \mathbf{Q}) \quad \forall \mathbb{Q} \in \mathcal{Q}$$ is infeasible in (\ref{prescription problem}). Furthermore, $(\hat{f}, \hat{x})$ is \textit{strongly optimal with respect to a nonempty set $\widetilde{\mathcal{Z}} \subseteq \mathcal{Q}$} (\textit{S}-optimal w.r.t. $\widetilde{\mathcal{Z}}$), if it is feasible in (\ref{prediction problem}) and any other \textit{feasible} pair $(\hat{f}', \hat{x}')$ satisfies 
$$\hat{f}'(\hat{x}'(\mathbf{Q}), \mathbf{Q}) \geq \hat{f}(\hat{x}'(\mathbf{Q}), \mathbf{Q}) \quad \forall \mathbb{Q} \in \widetilde{\mathcal{Z}} $$
\vspace{-9.5mm}\flushright$\square$ 
\end{definition}
	
Analogously, if $\widetilde{\mathcal{Z}} = \mathcal{Q}$, then the outlined prediction-prescription pair $(\hat{f}, \hat{x})$ is strongly optimal in a common sense. 
We show that (\ref{distributionally robust prescriptor}) is strongly optimal with respect to some set $\widetilde{\mathcal{Z}} \subseteq \mathcal{Q}$ for the prescription problem (\ref{prescription problem}). The following result holds. 

\begin{theorem} \label{theorem 4}
	Assume that $r > 0$ is a required exponential decay rate and Assumption \textbf{A1'} holds. Also, let
	\begin{equation} \label{eq: set Z 3}
	\widetilde{\mathcal{Z}} =\Big\{\mathbb{Q} \in \mathcal{Q}: q_{a,i} = q_{b,i} > 0 \quad \forall i \in \{1, \ldots, d_a\}, \; \forall a,b \in \mathcal{A}, \; a \neq b
	\Big\},
	\end{equation}
	where $q_{a,i}$ are rational and denote the probability that $c_a$ equals $z_{a,i}$. 	
	Then the prediction-prescription pair $(\hat{f}_{DR}, \hat{x}_{DR})$ with parameters $r_a$, $a \in \mathcal{A}$, given by equation (\ref{eq: radius}) with $\delta_a = \frac{1}{n}$ is \textit{S}-optimal w.r.t. $\widetilde{\mathcal{Z}}$ for (\ref{prescription problem}).
\begin{proof}
The proof, in a sense, reiterates the proof of Theorem \ref{theorem 2} with some minor changes.
First, we assume that there exists a prediction-prescription pair $(\hat{f}', \hat{x}')$, which is less conservative than~$(\hat{f}_{DR}, \hat{x}_{DR})$ at $\widetilde{\mathbf{Q}} = (\widetilde{\mathbb{Q}}_1, \ldots, \widetilde{\mathbb{Q}}_n)$ defined as in the proof of Theorem \ref{theorem 2}, i.e.,
\begin{equation} \label{eq: new weak optimaility prescription proof 1} \nonumber
\varepsilon := \hat{f}_{DR}(\hat{x}_{DR}(\widetilde{\mathbf{Q}}), \widetilde{\mathbf{Q}}) - \hat{f}'(\hat{x}'(\widetilde{\mathbf{Q}}), \widetilde{\mathbf{Q}}) > 0
\end{equation}
In view of Definition \ref{def: weak optimality 2}, we need to demonstrate that the pair $(\hat{f}', \hat{x}')$ is infeasible in (\ref{prescription problem}).  

In order to simplify the notations we set $\widetilde{\mathbf{x}} := \hat{x}_{DR}(\widetilde{\mathbf{Q}})$ and $\mathbf{x}' := \hat{x}'(\widetilde{\mathbf{Q}})$, where $\widetilde{\mathbf{x}}, \mathbf{x}' \in X$. 
Then we consider the out-of-sample prescription disappointment (\ref{eq: out-of-sample prescription}) given by:
$$\Pr\Big\{f(\hat{x}'(\widehat{\mathbf{Q}}), \mathbb{Q}^*) > \hat{f}'(\hat{x}'(\widehat{\mathbf{Q}}), \widehat{\mathbf{Q}}) \Big\}$$
and revise (\ref{eq: new weak optimality proof 3}) as follows:
\begin{equation} \label{eq: new weak optimality prescription proof 2} \nonumber
\begin{gathered}
f(\mathbf{x}', \mathbb{Q}^*) = \sum_{a \in \mathcal{A}} \mathbb{E}_{\mathbb{Q}^*_a}\{c_a\} x'_a > \sum_{a \in \mathcal{A}} \Big(\max_{\mathbb{Q}_a \in \widehat{\mathcal{Q}}_a}\mathbb{E}_{\mathbb{Q}_a}\{c_a\} - \varepsilon' \Big)x'_a = \\
\hat{f}_{DR}(\mathbf{x}', \widetilde{\mathbf{Q}}) - \varepsilon'\sum_{a \in \mathcal{A}} x'_a \geq \hat{f}_{DR}(\widetilde{\mathbf{x}}, \widetilde{\mathbf{Q}}) - \varepsilon' \sum_{a \in \mathcal{A}} x'_a = \hat{f}^{'}(\mathbf{x}', \widetilde{\mathbf{Q}}) + \varepsilon - \frac{\varepsilon}{nu} \sum_{a \in \mathcal{A}} x'_a \geq \hat{f}^{'}(\mathbf{x}', \widetilde{\mathbf{Q}}),
\end{gathered} 
\end{equation} 
where the second inequality follows from the definition of $\widetilde{\mathbf{x}}$.
We conclude that
$$\Pr\Big\{f(\hat{x}'(\widehat{\mathbf{Q}}), \mathbb{Q}^*) > \hat{f}'(\hat{x}'(\widehat{\mathbf{Q}}), \widehat{\mathbf{Q}}) \Big\} \geq \Pr\Big\{\widehat{\mathbf{Q}} = \widetilde{\mathbf{Q}} \Big\}$$
and the result follows from the second step in the proof of Theorem \ref{theorem 2}.	
\end{proof}
\end{theorem}

Finally, we formulate the result, which holds for any $d_a \in \mathbb{Z}_{>0}$, $a \in \mathcal{A}$; recall Theorem \ref{theorem 3}. 

\begin{theorem} \label{theorem 5}
	Assume that $r > 0$ is a required exponential decay rate and Assumption \textbf{A1''} holds. Also, let
	\begin{equation} \label{eq: set Z 4}
	\widetilde{\mathcal{Z}} = \Big\{\mathbb{Q} \in \mathcal{Q}: q_{a,i} = \begin{cases} q_{b,i} \in \mathbb{R}_{>0} \mbox{, if } i \leq d_{min} \\
	0, \mbox{ otherwise} \end{cases} \quad \forall a,b \in \mathcal{A}, \; a \neq b
	\Big\},
	\end{equation}
	where $q_{a,i}$ are rational and denote the probability that $c_a$ equals $z_{a,i}$. 	
	 Then the distributionally robust prescription rule (\ref{distributionally robust prescriptor}) with parameters $r_a$, $a \in \mathcal{A}$, given by equation (\ref{eq: radius}) with $\delta_a = \frac{1}{n}$ is \textit{S}-optimal w.r.t. $\widetilde{\mathcal{Z}}$ for (\ref{prescription problem}).
	
	\begin{proof} The result follows from Theorems \ref{theorem 3} and \ref{theorem 4}.
	\end{proof}
\end{theorem}

Theorems \ref{theorem 2} and \ref{theorem 4} establish that the prediction and prescription rules (\ref{distributionally robust predictor}) and (\ref{distributionally robust prescriptor}) are optimal in the asymptotic sense whenever the empirical marginal distributions $\widehat{\mathbb{Q}}_{a,T_a}$, $a \in \mathcal{A}$, have affinely dependent support sets but the same probability mass functions; recall equation (\ref{eq: set Z}). Admittedly, in practice this modeling assumption cannot be realized as the decision-maker is not able to control the initial data~set. Nevertheless, in our numerical experiments we attempt to relax the outlined ideal situation by leveraging a complete data set (\ref{data set}), i.e., with $T_a = T$, $a \in \mathcal{A}$, obtained from independent marginal distributions of the same type. If the parameter $T$ is sufficiently large, then the obtained empirical marginal distributions usually have a \textit{similar form} up to some affine transformation. We demonstrate numerically that in this case our approach outperforms several other benchmark approaches in terms of their out-of-sample performance; see Section \ref{subsec: results and discussion} and, in particular, Figure~\ref{fig: binomial distribution, Tmax} for further~details.

As we outlined before, most of the previous theoretical results focus on asymptotic performance guarantees. Alternatively, in the next section we propose potential improvements to the finite sample guarantee (\ref{eq: finite sample guarantee}) provided via the proof of Theorem \ref{theorem 1} and used in our numerical experiments.

\subsection{Improved finite sample guarantees} \label{subsec: improved finite sample}
In this section we assume that Assumption \textbf{A1'} holds and, thus, the conditions of Theorems \ref{theorem 2} and \ref{theorem 4} are satisfied. 
Following the proof of Theorem~\ref{theorem 1} we observe that the out-of-sample prediction disappointment for (\ref{distributionally robust predictor}) is bounded from above as follows:
\begin{equation} \nonumber
\Pr\Big\{f(\mathbf{x}, \mathbb{Q}^*) > \hat{f}_{DR}(\mathbf{x}, \widehat{\mathbf{Q}}) \Big \} \leq \sum_{a \in \mathcal{A}_{\mathbf{x}}}\mbox{\upshape Pr} \Big\{D_{KL}(\mathbb{Q}_a \; \Vert \; \mathbb{Q}^*_a) > r_a \Big\} \leq \sum_{a \in \mathcal{A}}(T_a + 1)^{d_a}e^{-T_a r_a},
\end{equation}
where the last inequality exploits the upper bound (\ref{eq: strong LDP 1}) and the fact that $\mathcal{A}_{\mathbf{x}} \subseteq \mathcal{A}$ for any $\mathbf{x} \in X$. 
However it can be argued that the upper bound (\ref{eq: strong LDP 1}) is not tight. 

Recently 
Mardia~et~al.~\cite{Mardia2020} has proposed some tighter upper bound that can be used to provide less conservative estimates of the parameters $r_a$, $a \in \mathcal{A}$. 
That is, by Theorem 3 in \cite{Mardia2020} for any $d_a \geq 2$ and $T_a \geq 2$ the following upper bound holds:
\begin{equation} \label{eq: strong LDP 2} \tag{$\mathcal{LDP}_2$}
\mbox{\upshape Pr} \Big\{D_{KL}(\mathbb{Q}_a \; \Vert \; \mathbb{Q}^*_a) > r_a \Big\} \leq \Big(\frac{3u_1}{u_2} \sum_{j = 0}^{d_a-2} K_{j-1} (\frac{e\sqrt{T_a}}{2\pi})^j\Big) e^{-T_a r_a},
\end{equation}
where $u_0 = \pi$, $u_1 = 2$, $K_{-1} = 1$,
\begin{equation} \nonumber
u_i = \begin{cases}\pi \times \frac{1 \times 3 \times \ldots \times (i-1)}{2 \times 4 \times \ldots \times i}, \mbox{ if } i \mbox{ is even and } i \geq 2 \\
2 \times \frac{2 \times 4 \times \ldots \times (i-1)}{1 \times 3 \times \ldots \times i}, \mbox{ if } i \mbox{ is odd and } i \geq 3
\end{cases} \mbox{ and } K_{j} = \prod_{i = 0}^{j}u_i \quad \forall j \geq 1
\end{equation}

In particular, it is illustrated in \cite{Mardia2020} that for most parameter settings the upper bound (\ref{eq: strong LDP 2}) is tighter than (\ref{eq: strong LDP 1}) and some other existing upper bounds; see, e.g., \cite{Agrawal2020}. By using the proof of Theorem \ref{theorem 1}, in our numerical experiments we select the parameters $r_a$, $a \in \mathcal{A}$, by setting the right-hand side of (\ref{eq: strong LDP 2}) equal to $\delta_a e^{-rT_{min}}$. In particular, if Assumption \textbf{A1'} holds and $\delta_a = \frac{1}{n}$, $a \in \mathcal{A}$, then it is rather easy to check that 
$$\Pr\Big\{f(\mathbf{x}, \mathbb{Q}^*) > \hat{f}_{DR}(\mathbf{x}, \widehat{\mathbf{Q}}) \Big \} \leq e^{-rT_{min}}$$
and both Theorems \ref{theorem 2} and \ref{theorem 4} remain valid with the aforementioned choice of $r_a$, $a \in \mathcal{A}$. 


\section{Computational study} \label{sec: comp study}

The primarily goal of this section is to explore the quality of distributionally robust decisions induced by the prescription rule (\ref{distributionally robust prescriptor}). In particular, we provide a detailed numerical comparison of our approach with Hoeffding bounds and some modification of the DRO approach described in \cite{VanParys2020}. Following the research question \textbf{Q3} we examine different parameter settings with respect to both the form of the nominal distribution and the sample size. 


For a given nominal distribution $\mathbb{Q}^* \in \mathcal{Q}$ and a vector of empirical marginal distributions $\widehat{\mathbf{Q}}$, we measure the quality of a prescription rule $\hat{x}$ by using a \textit{nominal relative loss}, that is,
\begin{eqnarray} \label{eq: nominal relative loss}
\rho(\mathbb{Q}^*, \widehat{\mathbf{Q}}) = \frac{f(\hat{x}(\widehat{\mathbf{Q}}), \mathbb{Q}^*)}{\min_{\mathbf{x} \in X} f(\mathbf{x}, \mathbb{Q}^*)}
\end{eqnarray}
In fact, the nominal relative loss (\ref{eq: nominal relative loss}) evaluates the out-of-sample performance of $\hat{x}(\widehat{\mathbf{Q}}) \in X$ under the nominal distribution~$\mathbb{Q}^*$. Ideally, we have $\rho = 1$, while in general $\rho \geq 1$ as long as the decision-maker operates with a limited information about $\mathbb{Q}^*$.

With respect to the nominal problem in (\ref{stochastic programming problem}), we consider two types of combinatorial optimization problems, namely, the shortest path and unweighted knapsack problems. 
As outlined in Proposition~\ref{proposition 1}, computation of (\ref{distributionally robust prescriptor}) requires solving $n$ univariate convex optimization problems and a unique deterministic MIP problem. For this reason, we do not consider some more complicated MIP problems, but instead focus on some qualitative insights implied by our construction of the nominal distribution and the training data set. In particular, some practical inference from estimation of the worst-case expected costs for each component of the cost vector $\mathbf{c}$ is provided.

The remainder of this section is organized as follows. In Section \ref{subsec: benchmark approaches} we design several benchmark approaches based on standard measure concentration inequalities or truncation of the training data~set. Section \ref{subsec: test instances} provides the test instances and parameter settings. Finally, in Section \ref{subsec: results and discussion} we report and discuss our numerical results. 

\subsection{Benchmark approaches} \label{subsec: benchmark approaches}


\textbf{Hoeffding bounds.} \looseness-1 
Following Example \ref{example 1} we define prediction and prescription rules as:
\begin{subequations}
	\begin{align}
	& \hat{f}_{hoef}(\mathbf{x}, \widehat{\mathbf{Q}}) = \sum_{a \in \mathcal{A}} \min\Big \{\frac{1}{T_a}\sum_{j = 1}^{T_a} \widehat{c}_{a,j} + \varepsilon_a; \overline{z}_a \Big \} x_a \label{hoeffding predictor} \\
	& \hat{x}_{hoef}(\widehat{\mathbf{Q}}) \in \argmin_{\mathbf{x} \in X} \hat{f}_{hoef}(\mathbf{x}, \widehat{\mathbf{Q}}),
	\label{hoeffding prescriptor}
	\end{align}
\end{subequations}

where $\varepsilon_a$, $a \in \mathcal{A}$, are given by equation (\ref{eq: eps hoeffding}) with $\delta_a = \frac{1}{n}$.

\textbf{Truncated DRO methods (DRO$_1$ and DRO$_2$).} We consider two alternative approaches, namely, DRO$_1$ and DRO$_2$, based on a truncation of the training data set (\ref{data set}). In the truncated DRO$_1$ approach we exploit the first $T_{min}$ random observations of $c_a$ for each $a \in \mathcal{A}$ and apply the model of Van Parys~et~al. \cite{VanParys2020} assuming that the obtained complete data set is generated from some joint distribution~$\mathbb{Q}^*$. In this model the decision-maker aims at minimizing its worst-case expected loss with respect to all probability distributions within a relative entropy ball centered at the \textit{joint empirical distribution} of the data, say, $\hat{\mathbb{J}}(T_{min})$. Then the outlined prediction and prescription rules can be described by their dual formulations as follows (we omit some minor technical details for brevity; see Proposition 2 in \cite{VanParys2020}):
\begin{subequations}
	\begin{align}
	& \hat{f}_{trunc}(\mathbf{x}, \hat{\mathbb{J}}) = \min_{\beta}\Big\{\beta - e^{-\widetilde{r}}\prod_{i = 1}^{\widetilde{d}} (\beta - (\widehat{\mathbf{c}}^{(i)})^\top \mathbf{x})^{\widehat{j}_{i}}: \; \beta \geq \sum_{a \in \mathcal{A}} \overline{z}_a x_a\Big\} \label{truncated predictor} \\
	& \hat{x}_{trunc}(\hat{\mathbb{J}}) \in \argmin_{\mathbf{x} \in X} \hat{f}_{trunc}(\mathbf{x}, \hat{\mathbb{J}})
	\label{truncated prescriptor}
	\end{align}
\end{subequations}
Here, $\widetilde{d} = \prod_{a \in \mathcal{A}}d_a$ is a number of possible realizations of $\mathbf{c}$, $\widehat{j}_{i}$ is an empirical probability that $\mathbf{c} = \widehat{\mathbf{c}}^{(i)}$ for $i \in \{1, \ldots, \widetilde{d}\}$ and $\widetilde{r}$ is a radius of the relative entropy ball centered at $\widehat{\mathbb{J}}$. Since the prediction rule (\ref{truncated predictor}) obeys finite sample guarantees similar to those obtained in Theorem\ref{theorem 1} (see Theorem 5 in~\cite{VanParys2020}), we compute $\widetilde{r}$ for some fixed decay rate $r > 0$ using the upper bound (\ref{eq: strong LDP 2}) adapted to the joint distributions. In particular, as the value of $\widetilde{d}$ is typically large, the sum in (\ref{eq: strong LDP 2}) is computed by eliminating zero terms in the sense of a floating-point precision.

\looseness-1 On the other hand, in the truncated DRO$_2$ approach we apply our prediction and prescriptions rules, (\ref{distributionally robust predictor}) and (\ref{distributionally robust prescriptor}), to the truncated data set. In other words, we set $T_a = T_{min}$ for each $a \in \mathcal{A}$. 

While the truncated DRO$_2$ approach preserves the complexity of the distributionally robust optimization problem (\ref{DRO problem}), the application of DRO$_1$ approach in our problem setting is substantially limited; recall our discussion in Section \ref{subsec: approach and contribution}.
In fact, the objective function in (\ref{truncated prescriptor}) is non-linear and the set of feasible decisions, $X$, is generally non-convex. Hence, computation of (\ref{truncated prescriptor}) results in a non-linear MIP problem. In our computational experiments we simply enumerate all decisions $\mathbf{x} \in X$ and solve the resulting univariate convex optimization problems in the parameter~$\beta$. In view of the discussion above, a further implementation of the truncated DRO$_1$ approach is restricted to comparatively small instances of the test problems.

\subsection{Test instances} \label{subsec: test instances}
\begin{figure}
	\centering
	\begin{tikzpicture}[scale=0.75,transform shape]
	\Vertex[x=0,y=0]{1}
	\Vertex[x=3,y=3]{2}
	\Vertex[x=3,y=0]{3}
	\Vertex[x=3,y=-3]{4}
	\Vertex[x=6,y=3]{5}
	\Vertex[x=6,y=0]{6}
	\Vertex[x=6,y=-3]{7}
	\Vertex[x=9,y=3]{8}
	\Vertex[x=9,y=0]{9}
	\Vertex[x=9,y=-3]{10}
	\Vertex[x=12,y=0]{11}
	\tikzstyle{LabelStyle}=[fill=white,sloped]
	\tikzstyle{EdgeStyle}=[post]
	\Edge(1)(2)
	\Edge(1)(3)
	\Edge(1)(4)
	\Edge(8)(11)
	\Edge(9)(11)
	\Edge(10)(11)
	\Edge(2)(5)
	\Edge(2)(6)
	\Edge(2)(7)
	\Edge(3)(5)
	\Edge(3)(6)
	\Edge(3)(7)
	\Edge(4)(5)
	\Edge(4)(6)
	\Edge(4)(7)
	\Edge(5)(8)
	\Edge(5)(9)
	\Edge(5)(10)
	\Edge(6)(8)
	\Edge(6)(9)
	\Edge(6)(10)
	\Edge(7)(8)
	\Edge(7)(9)
	\Edge(7)(10)
	\end{tikzpicture}
	\caption{\footnotesize A fully-connected layered graph with $h = 3$ intermediate layers and $w = 3$ nodes at each layer.}
	\label{fig: layered graph}
\end{figure}
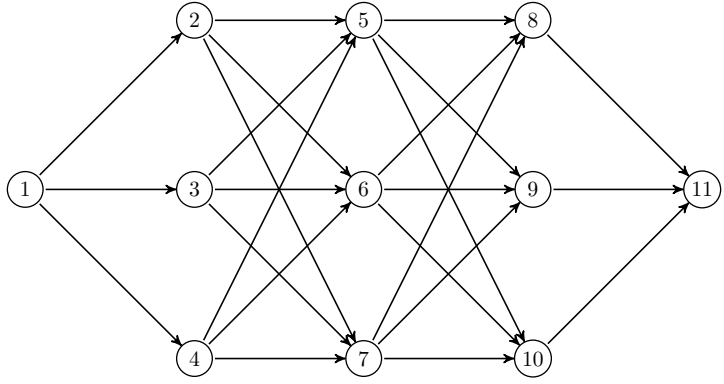

\textbf{Nominal problems.} As outlined earlier, we focus on two types of combinatorial optimization problems. First, we consider the shortest path problem in a fully-connected layered graph with $h$ intermediate layers and $w$ nodes at each layer $i \in \{1, \ldots, h\}$. The first and the last layers consist of unique nodes, which are the source and the destination nodes, respectively. An example with $h = w = 3$ is depicted in Figure \ref{fig: layered graph}. Specifically, we observe that the indexes $\mathcal{A} = \{1, \ldots, n\}$ are related to the set of directed arcs, whereas a decision $\mathbf{x} \in \{0, 1\}^n$ encodes a simple path between the source and destination nodes. Hence, $X$ is described by the standard path flow constraints \cite{Ahuja1988}, which explicit form is omitted for brevity.

Second, we consider a stochastic version of the unweighted knapsack problem (UNP), where the decision-maker attempts to minimize its expected loss by selecting a decision vector $\mathbf{x} \in \{0, 1\}^n$ with at least $K \in \{1, \ldots, n\}$ nonzero components. Formally, we consider the stochastic programming problem (\ref{stochastic programming problem}) with a set of feasible decisions given by:
\begin{equation} \label{feasible set: knapsack}
X_{UNP} = \Big\{\mathbf{x} \in \{0, 1\}^n: \sum_{i = 1}^n x_i \geq K \Big\}
\end{equation}
As a remark, this problem can be posed as an offline version of the best arm identification problem in online learning; see, e.g., \cite{Kalyanakrishnan2010} and recall our discussion of the application settings in Section \ref{subsec: approach and contribution}.

\looseness-1All experiments are performed on a PC with \textit{CPU i7-9700} and \textit{RAM 32 GB}. The deterministic versions of combinatorial optimization problems are solved in Python with \textit{CPLEX 12.10.0.0}. The dual formulations in (\ref{individual dual problems}) are solved using the function \textit{scipy.optimize.minimize}() and the method of Nelder-Mead with default parameters. In particular, we verify that the strong duality holds by solving the primal optimization problems (\ref{worst-case expectation problem}) with \textit{CVX 1.0.31}. 
 
\begin{table}
	
	\footnotesize
	\centering
	\onehalfspacing
	\begin{tabular}{c c c c c}
		\hline
		Joint distribution & Marginal distributions & Support & Mean & Variance\\
		\hline
		\multirow{2}{*}{Product of marginals} & \multirow{2}{*}{Binomial($p_a$, $d - 1$)} & \multirow{2}{*}{$c_a \in \{1, \ldots, d\}$} & \multirow{2}{*}{$(d - 1)p_a + 1$} & \multirow{2}{*}{$dp_a(1 - p_a)$} \\
		& & & & \\\hline
		\multirow{2}{*}{Multinomial($\mathbf{p}$, $\mathbf{d - 1}$)} & Binomial($p_a$, $d - 1$) & $c_a \in \{1, \ldots, d\}$, & \multirow{2}{*}{$(d - 1)p_a + 1$} & \multirow{2}{*}{$dp_a(1 - p_a)$}\\
		& with $\sum_{a \in \mathcal{A}} p_a = 1$ & $\sum_{a \in \mathcal{A}} c_a = d - 1 + |\mathcal{A}|$ & & \\\hline
	  \multirow{2}{*}{Product of marginals} & Discretization of normal & \multirow{2}{*}{$c_a \in \{1, \ldots, d\}$} & \multirow{2}{*}{$\sum_{i = 1}^d i q^*_{a, i}$} & $\sum_{i = 1}^d i^2 q^*_{a, i} - $\\
	  & with mean $\mu_a$ and variance $\sigma^2_a$ & & & $\Big(\sum_{i = 1}^d i q^*_{a, i}\Big)^2$\\
		\hline
	\end{tabular}
	\caption{\scriptsize The table summarizes discrete distributions that are used to model $\mathbb{Q}^*$. The values of binomial and multinomial distributions are shifted to the right by one. The support, mean and variance are component-wise. Parameters $q^*_{a, i}$, $a \in \mathcal{A}$, of the normal distributions are computed using (\ref{eq: prob normal discretized}).}
	\label{tab: data-generating distributions}
\end{table}

\textbf{Data-generating distributions.} \looseness-1 We assume that $\mathcal{S}_a = \{1, \ldots, d\}$, $d \in \mathbb{Z}_{>0}$, for each $a \in \mathcal{A}$, i.e., both Assumptions \textbf{A1} and \textbf{A1'} are satisfied. Then several classes of discrete distributions are examined; see Table \ref{tab: data-generating distributions}. The binomial and multinomial distributions are standard, but shifted to the right by one. The additional support constraint 
$$\sum_{a \in \mathcal{A}} c_a = d - 1 + |\mathcal{A}|$$
for the multinomial distribution is \textit{not known} to the decision-maker due to Assumption \textbf{A1}. 

In addition, we note that the variance of the binomial and multinomial distirbutions for each $a \in \mathcal{A}$ depends on the mean and, thus, cannot be controlled directly; see the last column of Table \ref{tab: data-generating distributions}. In this regard, we consider a naive discretization of univariate normal distributions, in which the variance for each $a \in \mathcal{A}$ can be controlled by an external parameter $\sigma_a$. That is, for a given $\mu_a \in [1, d]$ and $\sigma_a \geq 0$ the nominal probabilities $q^*_{a, i}$, $i \in \{1, \ldots, d\}$, are computed as follows:
\begin{equation} \label{eq: prob normal discretized}
q^*_{a, i} = \frac{C}{\sigma_a \sqrt{2 \pi}}\int_{i - 0.5}^{i + 0.5} e^{-\frac{1}{2} \Big(\frac{t - \mu_a}{\sigma_a}\Big)^2}\mbox{d}t,
\end{equation}
Specifically, $C > 0$ is a normalization constant and the integrals are computed numerically using the function \textit{scipy.integrate.quad()} in Python.

\textbf{Sample size.} In order to pick the values of $T_a$, $a \in \mathcal{A}$, we fix some $\widetilde{T}_{min} \in \mathbb{Z}_{>0}$ and set
\begin{equation} \nonumber
\widetilde{T}_{max} = \widetilde{T}_{min} + \Delta,
\end{equation}
where $\Delta \in \mathbb{Z}_{>0}$ is some positive integer constant. By leveraging the values of $\widetilde{T}_{min}$ and $\widetilde{T}_{max}$ we set the values of $T_a$, $a \in \mathcal{A}$, in three different ways; see Table \ref{tab: choice of T}.

 The intuition behind the choice of $T_a$ is as follows. The Binomial$_1$ distribution sorts the mean values of $T_a$, $a \in \mathcal{A}$, in an increasing order with respect to their nominal expected costs. Oppositely, the Binomial$_2$ distribution sorts the mean values of $T_a$, $a \in \mathcal{A}$, in the decreasing order. In other words, in the former situation the components with higher expected costs can be observed sufficiently often; in the latter situation the same holds for the components with lower expected costs. Finally, the uniform distribution is somewhat in the middle between the Binomial$_1$ and Binomial$_2$ distributions. In the next section we demonstrate that the distribution of $T_a$, $a \in \mathcal{A}$, substantially affects the out-of-sample tests both for Hoeffding bounds and the baseline DRO approach induced by (\ref{distributionally robust prescriptor}).

\begin{table}	
	\footnotesize
	\centering
	\onehalfspacing
	\begin{tabular}{c c c c c}
		\hline
		Distribution of $T_a$, $a \in \mathcal{A}$ & Support & Parameters \\
		\hline
			\multirow{2}{*}{Uniform} & 	\multirow{2}{*}{$[\widetilde{T}_{min}, \widetilde{T}_{max}]$} & 	\multirow{2}{*}{-} \\
			& & \\\hline
		\multirow{2}{*}{Binomial$_1$($p_a, \widetilde{T}_{max} - \widetilde{T}_{min}$)} & \multirow{2}{*}{$[\widetilde{T}_{min}, \widetilde{T}_{max}]$} & \multirow{2}{*}{$p_a = \frac{\mathbb{E}_{\mathbb{Q}_a^*}\{c_a\} - \min_{b \in \mathcal{A}}\mathbb{E}_{\mathbb{Q}_b^*}\{c_b\}}{\max_{b \in \mathcal{A}}\mathbb{E}_{\mathbb{Q}_b^*}\{c_b\} - \min_{b \in \mathcal{A}}\mathbb{E}_{\mathbb{Q}_b^*}\{c_b\}}$} \\
		& & \\\hline
		\multirow{2}{*}{Binomial$_2$($p_a, \widetilde{T}_{max} - \widetilde{T}_{min}$)} & \multirow{2}{*}{$[\widetilde{T}_{min}, \widetilde{T}_{max}]$} & \multirow{2}{*}{$p_a = 1 - \frac{\mathbb{E}_{\mathbb{Q}_a^*}\{c_a\} - \min_{b \in \mathcal{A}}\mathbb{E}_{\mathbb{Q}_b^*}\{c_b\}}{\max_{b \in \mathcal{A}}\mathbb{E}_{\mathbb{Q}_b^*}\{c_b\} - \min_{b \in \mathcal{A}}\mathbb{E}_{\mathbb{Q}_b^*}\{c_b\}}$} \\
		& & \\\hline
	\end{tabular}
	\caption{\scriptsize The table summarizes the ways to pick $T_a$, $a \in \mathcal{A}$. The binomial distributions are shifted to the right by $\widetilde{T}_{min}$.}
	\label{tab: choice of T}
\end{table}

\subsection{Results and discussion} \label{subsec: results and discussion}
In view of the discussion above, we compare the aforementioned solution approaches in terms of the nominal relative loss (\ref{eq: nominal relative loss}). Specifically, we compute the average relative loss and median absolute deviations around the mean (MADs) over $N_0 = 200$ randomly generated test instances. 

Next, we fix a confidence level $\alpha = e^{-rT_{min}}$ instead of the exponential decay rate, $r$, for convenience and set $\alpha = 0.05$ in all experiments. 
\looseness-1 In addition, we set $d_a = d = 50$ and $\delta_a = \frac{1}{n}$ for each $a \in \mathcal{A}$; recall the conditions of Theorems \ref{theorem 2} and \ref{theorem 4}. The default parameters of the shortest path problem are given by $h = 7$ and $w = 4$, if other is not specified. Finally, we set $n = 100$ and customize the parameter $K$ in the unweighted knapsack problem. 

\begin{figure}
	\begin{subfigure}[b]{0.5\textwidth}
		\centering
		\captionsetup{justification=centering}
		\includegraphics[width = \linewidth]{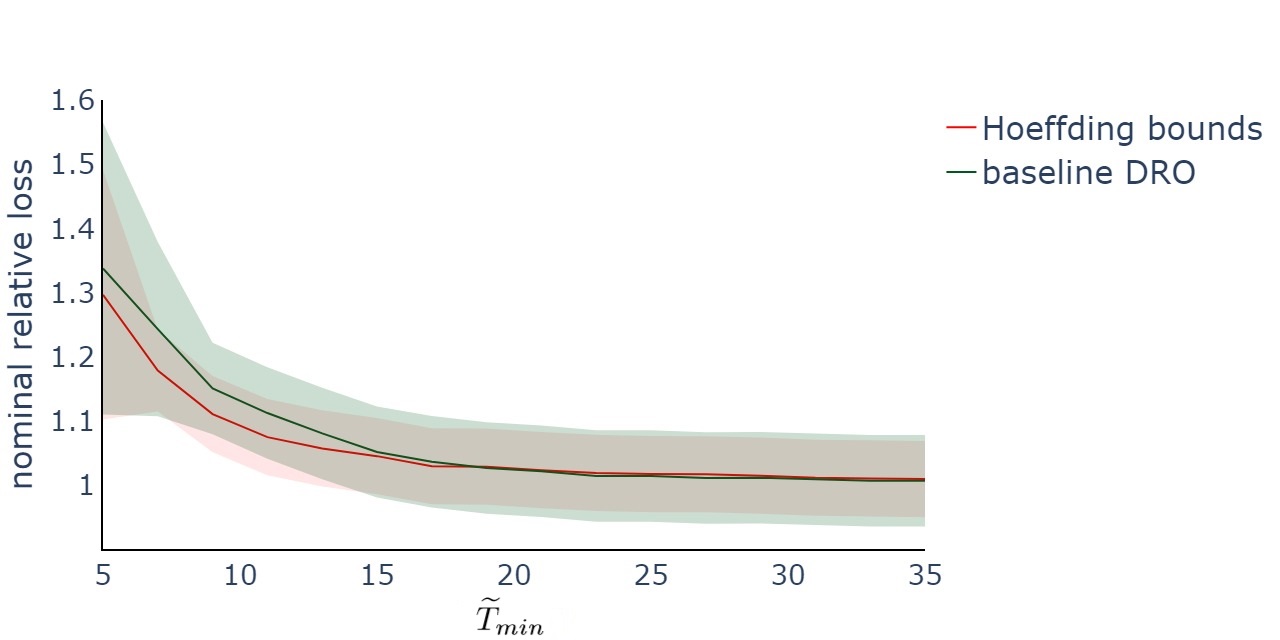}
		\caption{Binomial distribution.}	
		\label{fig: binomial, Tmin}
	\end{subfigure}
	\begin{subfigure}[b]{0.5\textwidth}
		\centering
		\captionsetup{justification=centering}	
		\includegraphics[width = \linewidth]{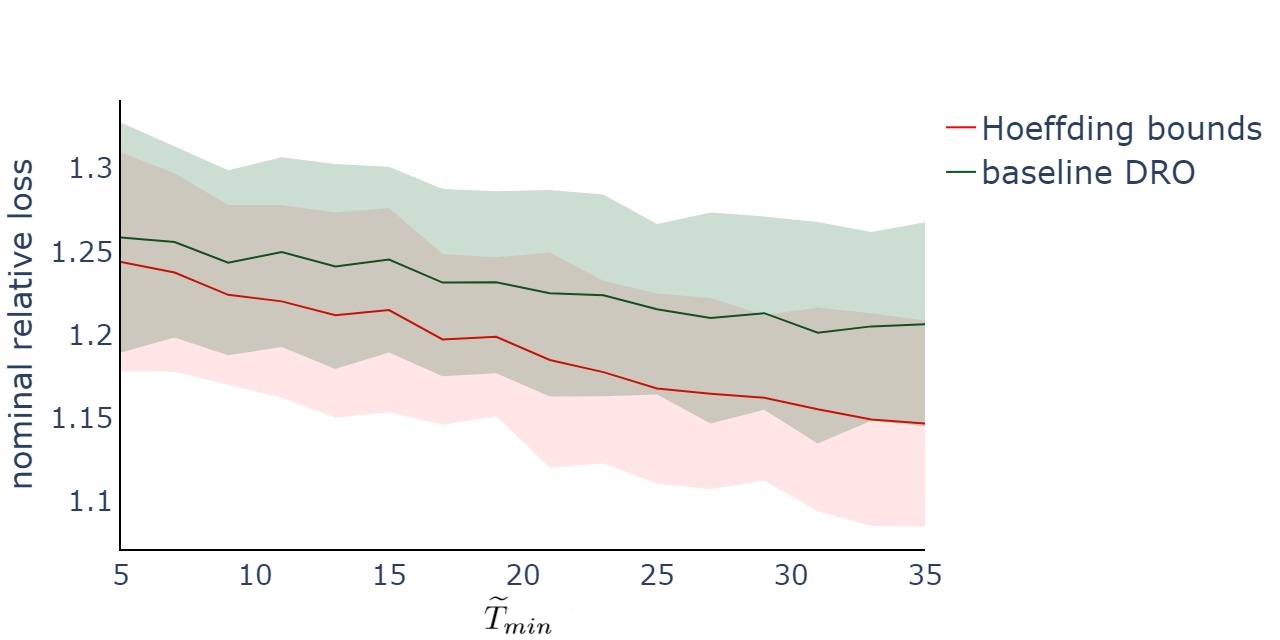}
		\caption{Multinomial distribution.}
		\label{fig: multinomial, Tmin}
	\end{subfigure}
	\caption{Average relative loss (\ref{eq: nominal relative loss}) and MADs as a function of $\widetilde{T}_{min}$, $\widetilde{T}_{min} \in \{5, 7 \ldots, 35\}$, with $\Delta = 10$ under the (a) binomial and (b) multinomial distributions. The distribution of $T_a$, $a \in \mathcal{A}$, is uniform, the parameters of the graph are given by $h = 7$ and $w = 4$.}
	\label{fig: binomial and multinomial, Tmin}
\end{figure}

\begin{figure}
	\begin{subfigure}[b]{0.5\textwidth}
		\centering
		\captionsetup{justification=centering}
		\includegraphics[width = \linewidth]{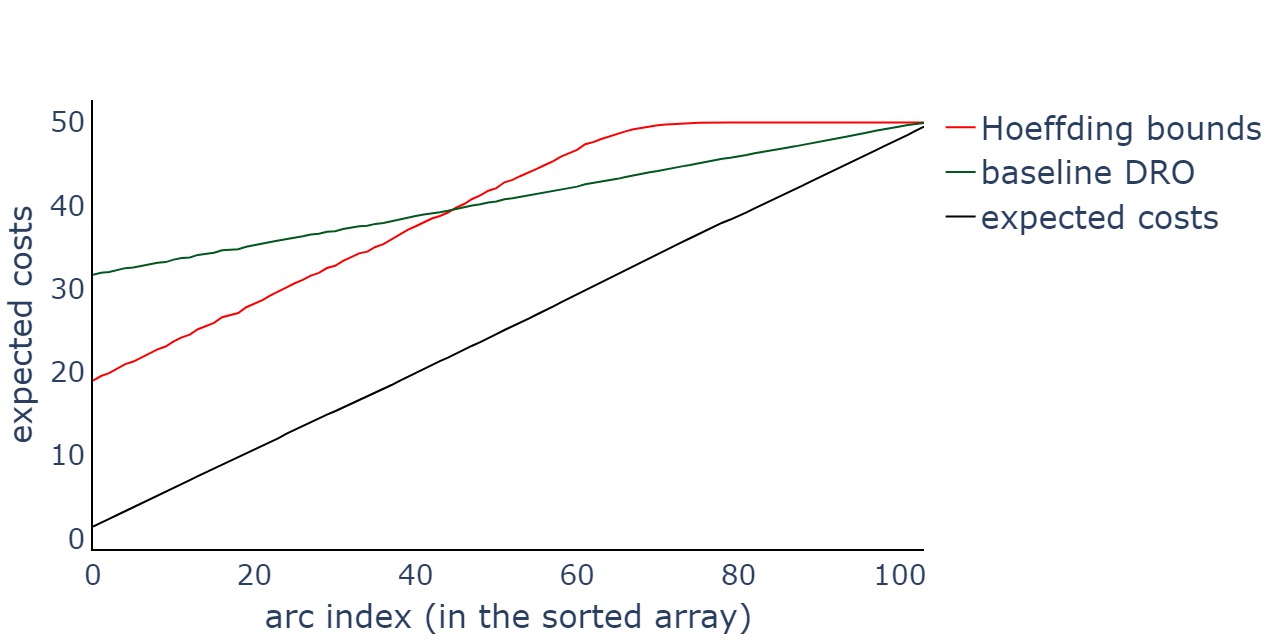}
		\caption{Binomial distribution.}	
		\label{fig: binomial, costs}
	\end{subfigure}
	\begin{subfigure}[b]{0.5\textwidth}
		\centering
		\captionsetup{justification=centering}	
		\includegraphics[width = \linewidth]{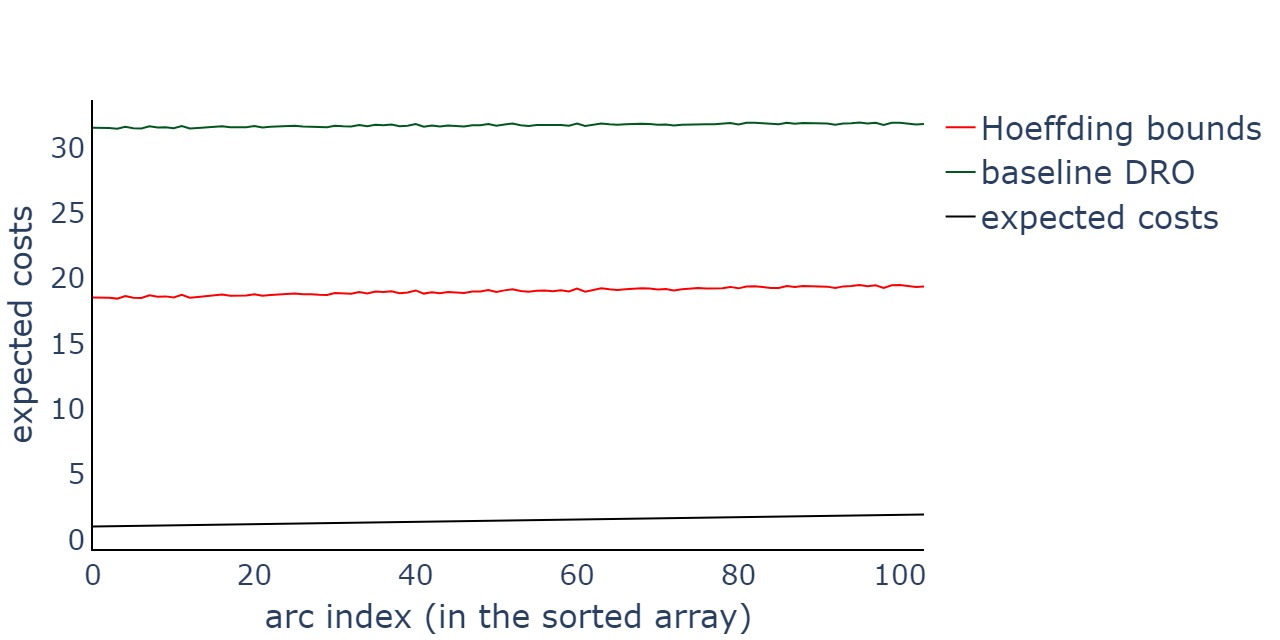}
		\caption{Multinomial distribution.}
		\label{fig: multinomial, costs}
	\end{subfigure}
	\caption{Average nominal expected costs (in an increasing order) and the associated worst-case expected costs under the (a) binomial and (b) multinomial distributions with $\widetilde{T}_{min} = 25$, $\Delta = 10$. The distribution of $T_a$, $a \in \mathcal{A}$, is uniform, the parameters of the graph are given by $h = 7$ and $w = 4$.}
	\label{fig: binomial and multinomial, costs}
\end{figure}

\subsubsection{Results for the shortest path problem}
\textbf{Dependence on the form of the nominal distribution.} 
In this paragraph we assume that the distribution of $T_a$, $a \in \mathcal{A}$, is uniform and compare the baseline DRO approach (\ref{distributionally robust predictor})-(\ref{distributionally robust prescriptor}) with Hoeffding bounds (\ref{hoeffding predictor})-(\ref{hoeffding prescriptor}). First, we consider the nominal relative loss~(\ref{eq: nominal relative loss}) as a function of $\widetilde{T}_{min}$ with $\Delta = 10$ under the binomial and multinomial distributions; see Figures \ref{fig: binomial, Tmin} and \ref{fig: multinomial, Tmin}, respectively. In particular, the parameters $p_a$ of binomial (multinomial) distributions for each $a \in \mathcal{A}$ are uniformly distributed over the interval $[0, 1]$; for the multinomial distribution the sum of $p_a$, $a \in \mathcal{A}$, is additionally normalized to one. 

As expected, the out-of-sample performance decreases as a function of $\widetilde{T}_{min}$; recall that the confidence level $\alpha$ is fixed. Furthermore, both methods provide similar results for the binomial distribution, but Hoeffding bounds demonstrate better out-of-sample performance in the case of the multinomial distribution. Some intuition behind this fact can be provided in terms of the worst-case expected costs. That is, in Figures \ref{fig: binomial, costs} and \ref{fig: multinomial, costs} we consider the nominal expected costs (sorted in an increasing order) and the associated worst-case expected costs averaged over $N_0 = 200$ instances with $\widetilde{T}_{min} = 25$ and $\Delta = 10$ . The key observation is that the lower expected costs are better estimated using Hoeffding bounds, while the higher costs are better estimated using the baseline DRO approach. As outlined in Table \ref{tab: data-generating distributions}, the multinomial distribution has an additional support constraint, which implies lower expected costs than those for the binomial distribution. This observation provides a practical evidence behind the plots in Figures \ref{fig: binomial, Tmin} and \ref{fig: multinomial, Tmin}.

In the second experiment we set $\mu_a$ uniformly distributed over $[1, d]$ and $\sigma_a = \sigma$ for each $a \in \mathcal{A}$; recall Table \ref{tab: data-generating distributions}. We consider the nominal relative loss (\ref{eq: nominal relative loss}) as a function of $\sigma$ under the discretized normal distribution with $\widetilde{T}_{min} = 25$ and $\Delta = 10$; see Figure \ref{fig: normal distribution, variance}. 
Naturally, the quality of obtained solutions decreases as the variance increases. However the baseline DRO approach and Hoeffding bounds demonstrate similar out-of-sample performance across the considered values of variance. Therefore, no general conclusions are made regarding a comparison of the considered solution approaches. 

\begin{figure}
	\begin{center}
		\includegraphics[width = 0.5\linewidth]{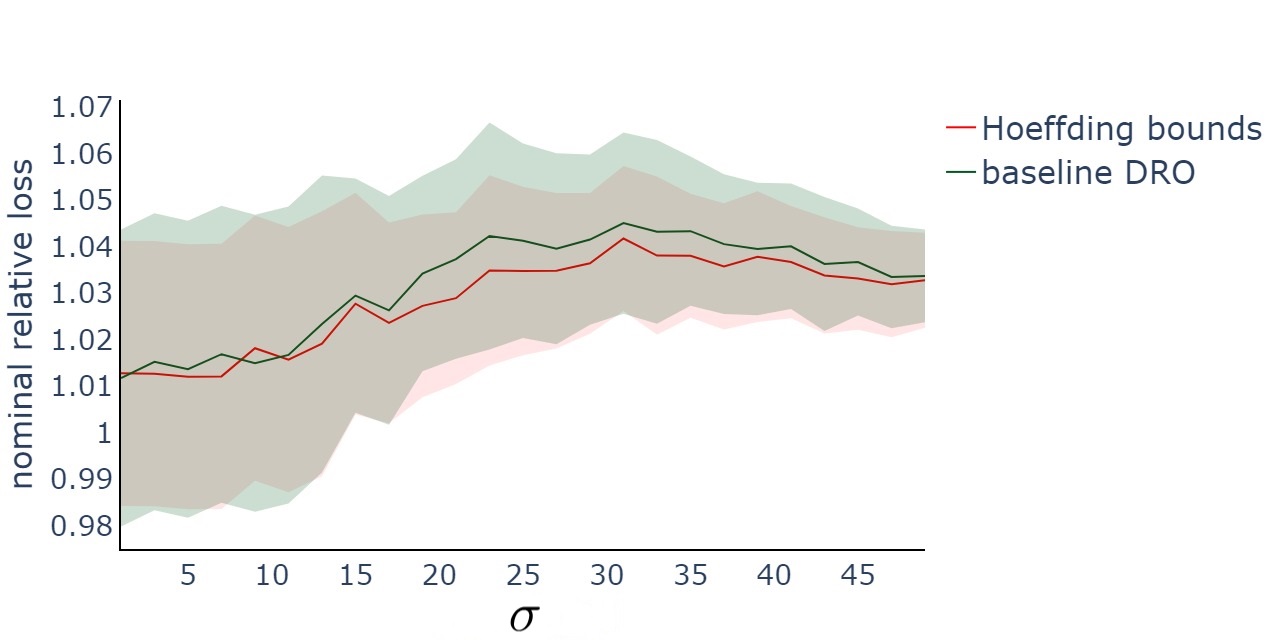}
		\caption{Average relative loss (\ref{eq: nominal relative loss}) and MADs as a function of $\sigma$, $\sigma \in \{1,3, \ldots, 49\}$, under the discretized normal distribution with $\widetilde{T}_{min} = 25$ and $\Delta = 10$ . The distribution of $T_a$, $a \in \mathcal{A}$, is uniform, the parameters of the graph are given by $h = 7$ and $w = 4$.}
		\label{fig: normal distribution, variance}
	\end{center}
\end{figure}

\begin{figure}
	\begin{subfigure}[b]{0.5\textwidth}
		\centering
		\captionsetup{justification=centering}
		\includegraphics[width = \linewidth]{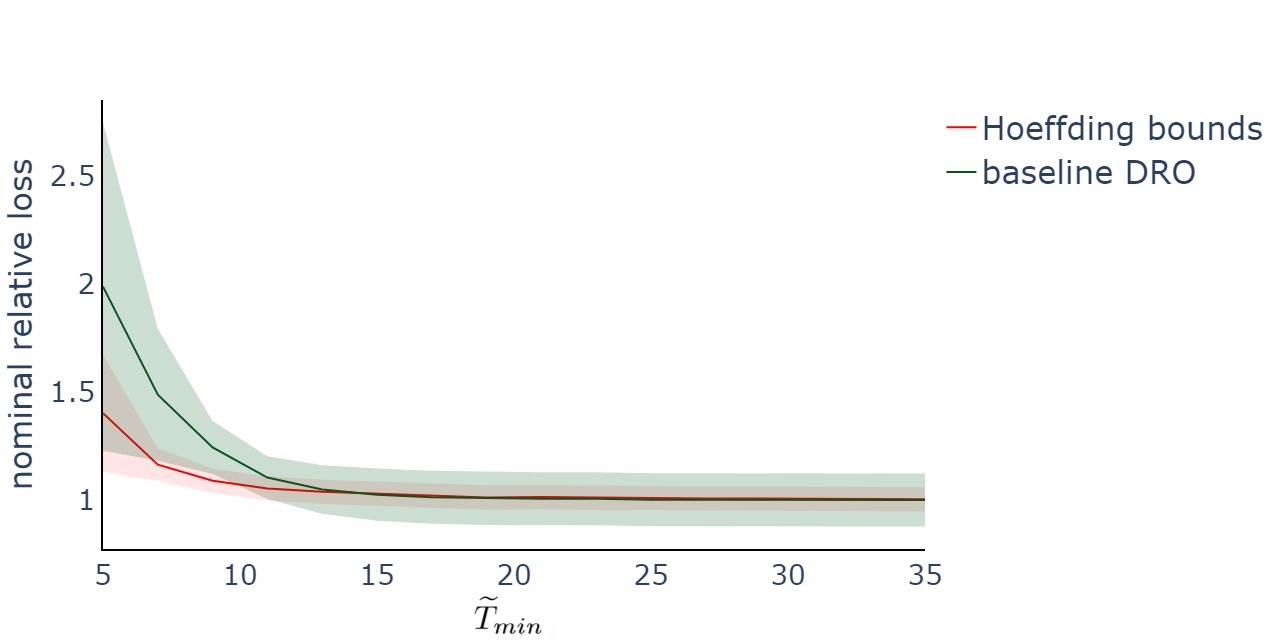}	
		\caption{Binomial distribution with Binomial$_1$ $T_a$, $a \in \mathcal{A}$.}	
		\label{fig: binomial T}
	\end{subfigure}
	\begin{subfigure}[b]{0.5\textwidth}
		\centering
		\captionsetup{justification=centering}
		\includegraphics[width = \linewidth]{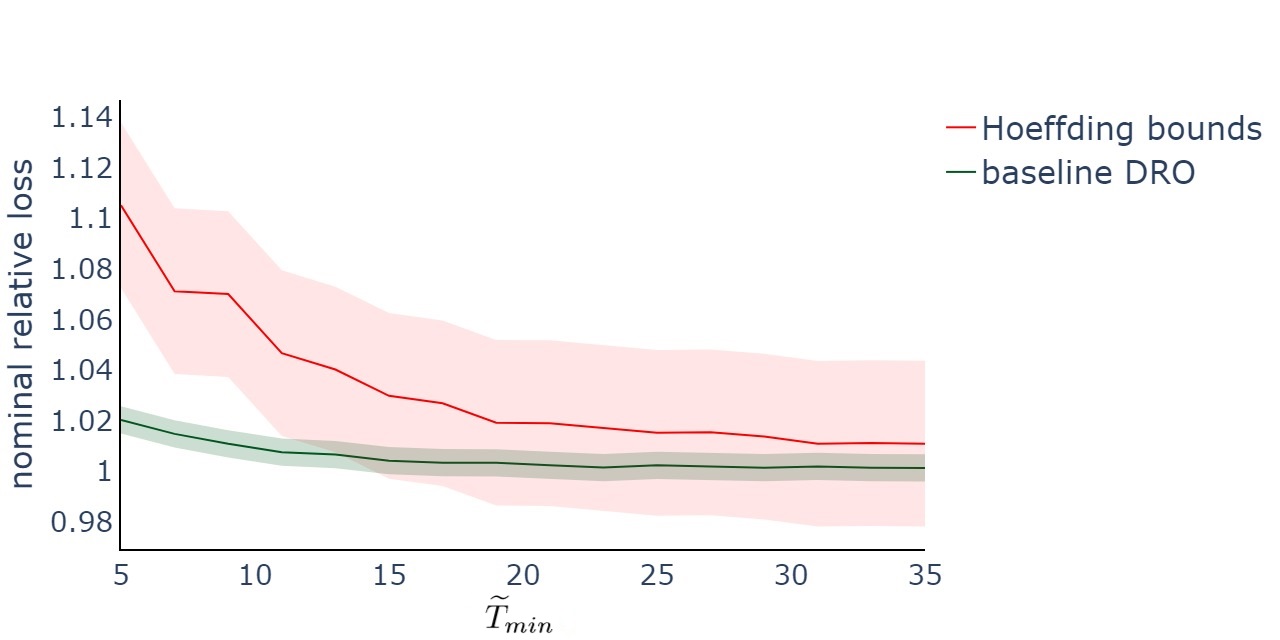}
		\caption{Binomial distribution with Binomial$_2$ $T_a$, $a \in \mathcal{A}$.}
		\label{fig: binomial T reversed}
	\end{subfigure}
	\caption{Average relative loss (\ref{eq: nominal relative loss}) and MADs as a function of $\widetilde{T}_{min}$, $\widetilde{T}_{min} \in \{5, 7 \ldots, 35\}$, with $\Delta = 10$ under the binomial distribution. The distributions of $T_a$, $a \in \mathcal{A}$, are assumed to be (a) Binomial$_1$ and (b) Binomial$_2$; the parameters of the graph are given by $h = 7$ and $w = 4$}
	\label{fig: binomial with binomial T}
\end{figure}

\textbf{Dependence on the sample size.} 
In this paragraph the data-generating distribution is supposed to be binomial with the parameters $p_a$, $a \in \mathcal{A}$, uniformly distributed over $[0, 1]$. We explore how the distribution of $T_a$, $a \in \mathcal{A}$, affects the out-of-sample performance of the baseline DRO approach and Hoeffding bounds. Specifically, we consider the nominal relative loss (\ref{eq: nominal relative loss}) as a function of $\widetilde{T}_{min}$ with $\Delta = 10$ under the Binomial$_1$ and Binomial$_2$ distributions from Table \ref{tab: choice of T}. The corresponding plots are depicted in Figures \ref{fig: binomial T} and \ref{fig: binomial T reversed}, respectively.

We conclude that the DRO approach outperforms Hoeffding bounds only in the case of Binomial$_2$ distribution.
This observation can be motivated in the context of Figure \ref{fig: binomial, costs}. That is, under the Binomial$_2$ distribution the estimates of lower expected costs become less conservative and, thus, the worst-case expected costs obtained from the DRO approach are, in a sense, more consistent with the nominal expected costs. The same holds for Hoeffding bounds and higher expected costs under the Binomial$_1$ distribution.

As we outlined in Section \ref{subsec: approach and contribution}, the Binomial$_2$ distribution arises naturally, if the data is collected by trial and errors through multiple decision epochs. More precisely, if the decision-maker takes a decision $\mathbf{x} \in X 
\subseteq \{0, 1\}^n$ and observes only the costs associated with this decision, then lower expected costs are more preferable for the decision-maker due to the objective criterion in (\ref{stochastic programming problem}). This observation reveals some novel motivation behind the proposed distributionally robust optimization approach.

In the next experiment we assume that the distribution of $T_a$, $a \in \mathcal{A}$, is uniform and explore the nominal relative loss (\ref{eq: nominal relative loss}) as a function of $\Delta$ for $\widetilde{T}_{min} = 10$ ; see Figure \ref{fig: binomial distribution, Tmax}. 
We note that for $\Delta = 0$ the nominal relative loss (\ref{eq: nominal relative loss}) under the baseline DRO approach is close to one despite the number of samples is sufficiently small, i.e., we have $T_a = 10$ for each $a \in \mathcal{A}$. 

Some intuition behind this result can be provided by Theorems \ref{theorem 2} and \ref{theorem 4}. That is, if $\Delta = 0$, then all empirical marginal distributions contain the same number of samples, $\widetilde{T}_{min}$, and are obtained from nominal marginal distributions of the same type, i.e., binomial distributions. We expect that in this case the empirical marginal distributions have, in a sense, a similar form; recall our discussion after the proof of Theorem \ref{theorem 4} in Section \ref{subsec: prescriptor}. Therefore, in view of Theorems \ref{theorem 2} and \ref{theorem 4}, we may anticipate that the aforementioned parameters setting provides some advantage to the baseline DRO approach. 

\begin{figure}[htp]
	\begin{center}
		\includegraphics[width = 0.5\linewidth]{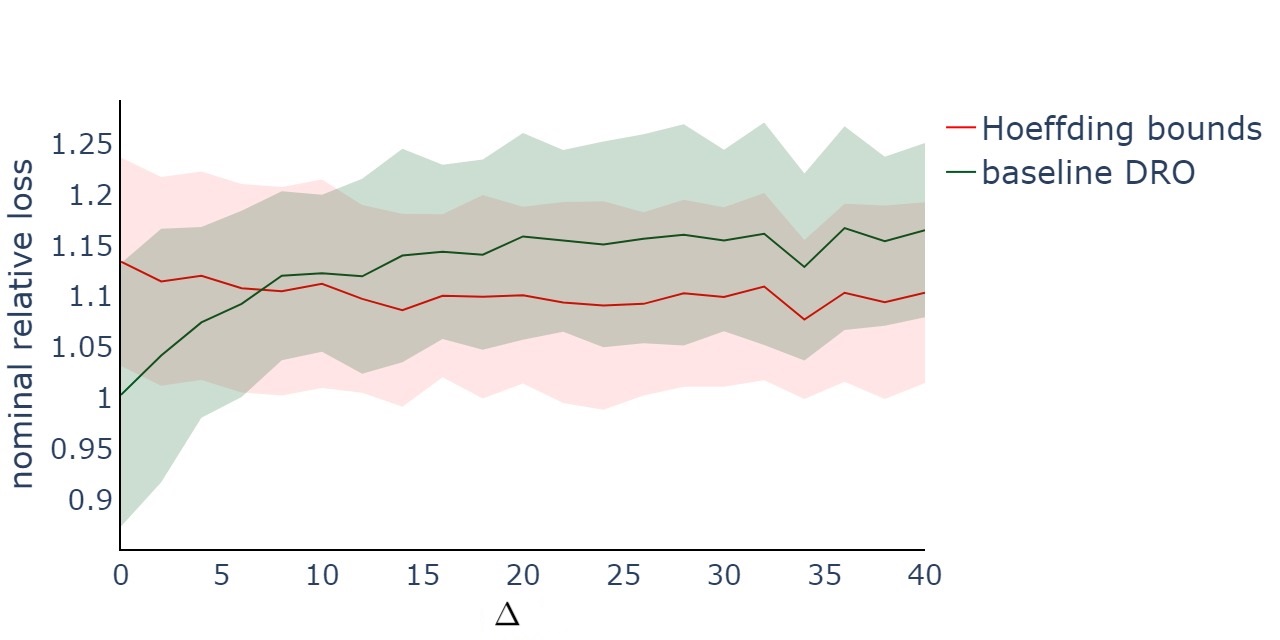}
		\caption{Average relative loss (\ref{eq: nominal relative loss}) and MADs as a function of $\Delta$, $\Delta \in \{0, 2, \ldots, 40\}$, under the binomial distribution with $\widetilde{T}_{min} = 10$. The distribution of $T_a$, $a \in \mathcal{A}$, is uniform, the parameters of the graph are given by $h = 7$ and $w = 4$. }
		\label{fig: binomial distribution, Tmax}	
	\end{center}
\end{figure}

Surprisingly, the increase of $\Delta$ and, thus, the use of additional data for particular components of $\mathbf{c}$ does not improve the out-of-sample performance of the baseline DRO approach; recall Figure \ref{fig: binomial distribution, Tmax}. In this regard, we make the following additional observations:
\begin{itemize}
	\item The increase of $\Delta$ implies a larger fluctuation of the worst-case expected costs obtained from both solution approaches.
	\item \looseness-1 In the case of the baseline DRO approach the worst-case expected costs increase rather \textit{slowly} in average as a function of the nominal expected costs; recall Figure \ref{fig: binomial, costs}. Hence, a large fluctuation of the worst-case expected costs may lead to non-uniform estimates of the nominal expected costs and, thus, to low-quality decisions in terms of the nominal relative loss.
	\item In the case of Hoeffding bounds the worst-case expected costs increase \textit{faster} as a function of the nominal costs. For this reason, the fluctuations in the case of sufficiently large $\Delta$ can be eliminated, in a sense, by less conservative estimates of the nominal expected costs; see Figure \ref{fig: binomial distribution, Tmax} and note that with the increase of $\Delta$ the sample size also tends to grow. 
	\item Furthermore, if $\Delta = 0$, then the worst-case expected costs obtained from Hoeffding bounds usually achieve the upper bound, $d$, which results in a poor out-of-sample performance of (\ref{hoeffding predictor})-(\ref{hoeffding prescriptor}). 
\end{itemize}

 We conclude that the baseline DRO approach outperforms Hoeffding bounds whenever the relative gap between sample sizes, $\Delta$, is sufficiently small. The obtained results for $\Delta = 0$ are extended in the next paragraph, where the truncated DRO$_1$ and DRO$_2$ methods are examined. 


\begin{figure}[htp]
	\begin{center}
		\includegraphics[width = 0.5\linewidth]{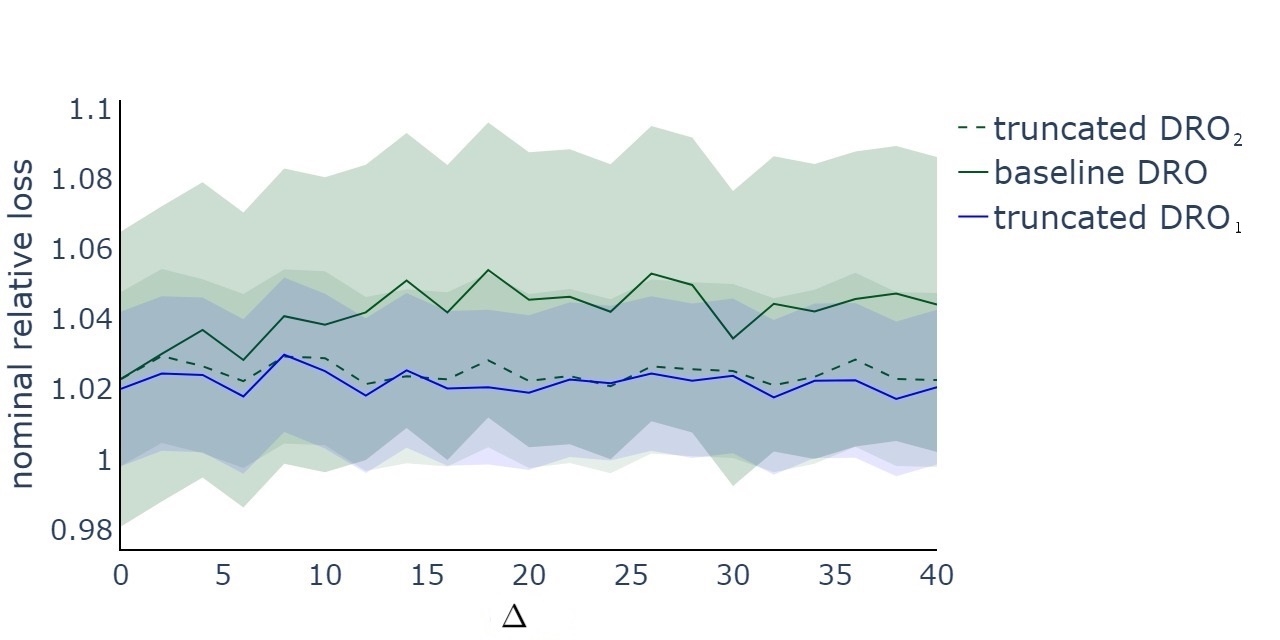}
		\caption{Average relative loss (\ref{eq: nominal relative loss}) and MADs as a function of $\Delta$, $\Delta \in \{0, 2, \ldots, 40\}$, under the discretized normal distribution with $\widetilde{T}_{min} = 10$ and $\sigma = \frac{d}{4}$. The distribution of $T_a$, $a \in \mathcal{A}$, is uniform and the parameters of the graph are given by $w = h = 3$.}
		\label{fig: normal distribution, sigma = d/4, T max}
	\end{center}
\end{figure}

\begin{figure}[htp]
	\begin{subfigure}[b]{0.5\textwidth}
		\centering
		\captionsetup{justification=centering}
		\includegraphics[width = \linewidth]{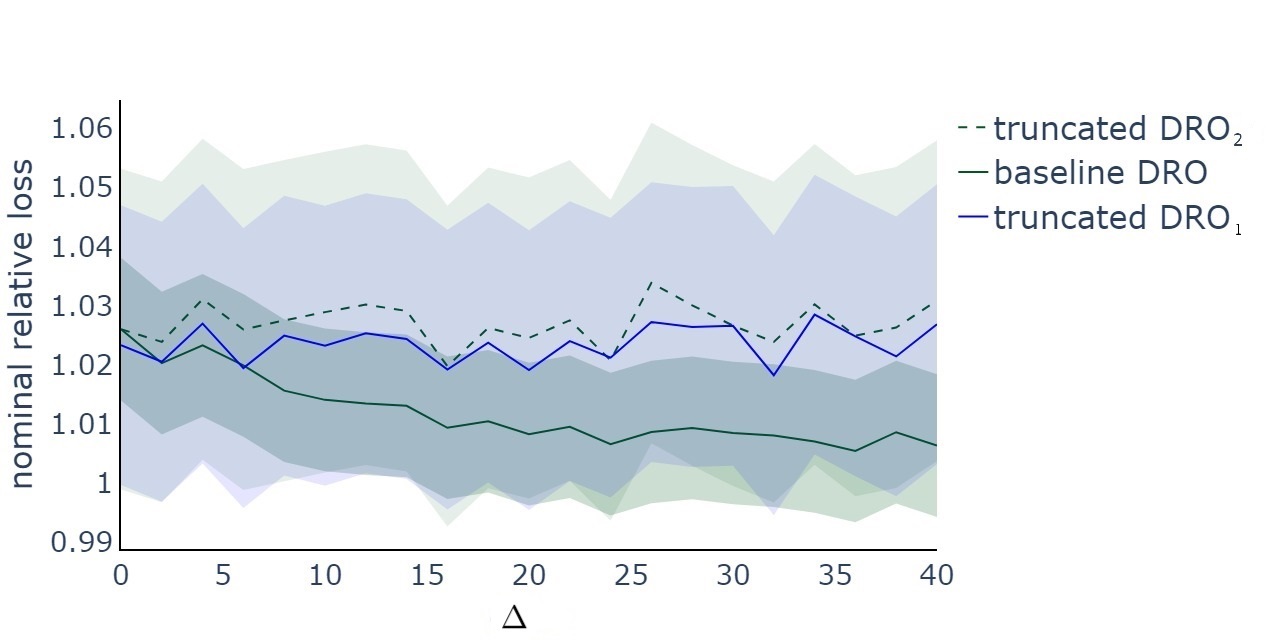}
		\caption{Normal distribution with Binomial$_2$ $T_a$, $a \in \mathcal{A}$, and $\sigma=\frac{d}{4}$.}
		\label{fig: normal, reversed T}
	\end{subfigure}
	\begin{subfigure}[b]{0.5\textwidth}
		\centering
		\captionsetup{justification=centering}
		\includegraphics[width = \linewidth]{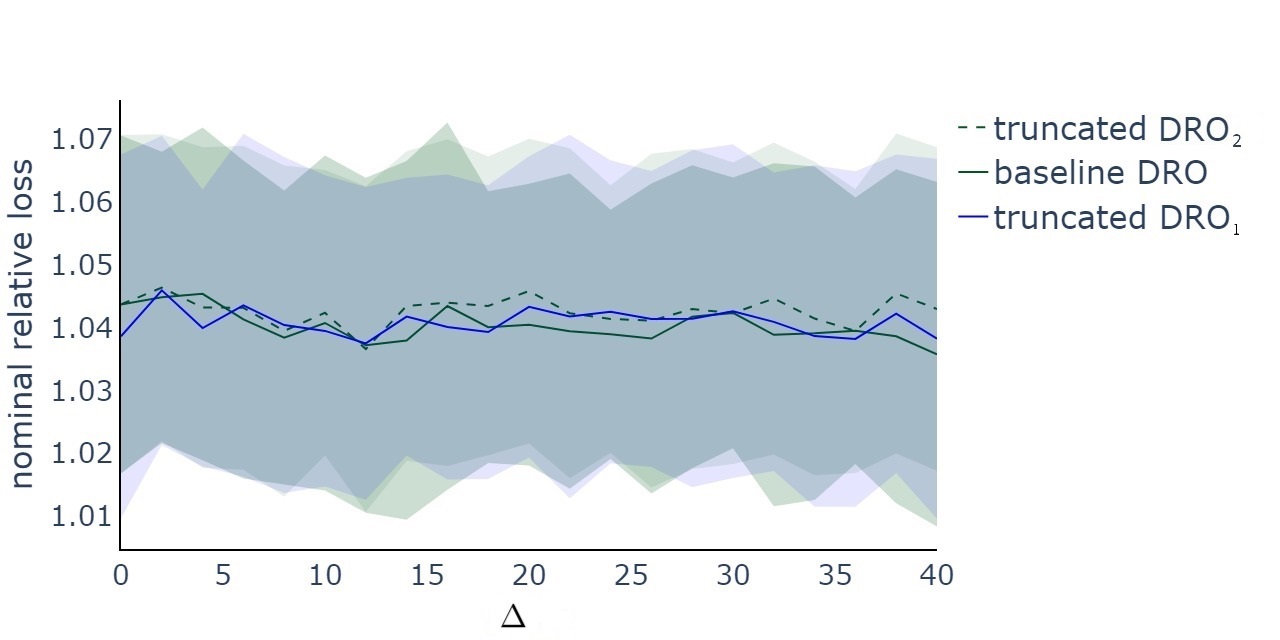}	
		\caption{Normal distribution with uniform $T_a$, $a \in \mathcal{A}$, and $\sigma=\frac{3d}{4}$.}	
		\label{fig: normal, sigma = 3d/4}
	\end{subfigure}
	\caption{Average relative loss (\ref{eq: nominal relative loss}) and MADs as a function of $\Delta$, $\Delta \in \{0, 2, \ldots, 40\}$, under the discretized normal distribution with (a) Binomial$_2$ $T_a$, $a \in \mathcal{A}$, $\sigma=\frac{d}{4}$ and (b) uniform $T_a$, $a \in \mathcal{A}$, $\sigma=\frac{3d}{4}$. In both cases we set $\widetilde{T}_{min} = 10$ and $w = h = 3$.}
	\label{fig: normal distribution, sigma = 3d/4 and reversed T}
\end{figure}

\textbf{Truncation of the data.}
\looseness-1 In this paragraph we assume that the nominal distribution is discretized normal with the parameters $\mu_a$ uniformly distributed over $[1, d]$ and $\sigma_a = \sigma$ for each $a \in \mathcal{A}$. Following our discussion of the truncated DRO$_1$ approach in Section \ref{subsec: benchmark approaches} we consider comparatively small instances of the shortest path problem with $h = w = 3$.

 First, we assume that the distribution of $T_a$, $a \in \mathcal{A}$, is uniform. The dependence on $\Delta$ is of interest since the truncated methods do not take into account a part of the data. That is, in Figure \ref{fig: normal distribution, sigma = d/4, T max} we depict the nominal relative loss (\ref{eq: nominal relative loss}) as a function of $\Delta$ with $\widetilde{T}_{min} = 10$ and $\sigma = \frac{d}{4}$. As a remark, the results for the truncated DRO$_1$ and DRO$_2$ approaches may also depend on $\Delta$ as it is not necessarily the case that $\widetilde{T}_{min} = T_{min}$; recall Table \ref{tab: choice of T}. 
	
The numerical results indicate that both methods with the truncated data set substantially outperform the baseline DRO approach for all considered values of parameter $\Delta$. This observation for the truncated DRO$_2$ approach is explained via the previous experiment; recall Figure \ref{fig: binomial distribution, Tmax}. At the same time, the truncated DRO$_1$ approach demonstrates a similar out-of-sample performance as the DRO$_2$ approach. Intuitively, the latter observation can be also provoked by the fact that the DRO$_1$ approach is asymptotically optimal; see Theorem 7 in \cite{VanParys2020}.

Finally, we note that the baseline DRO approach may outperform the DRO$_1$ and DRO$_2$ methods, e.g., if the parameters $T_a$, $a \in \mathcal{A}$, are governed by the Binomial$_2$ distribution or the variance $\sigma$ is relatively large; see Figures \ref{fig: normal, reversed T} and \ref{fig: normal, sigma = 3d/4}, respectively. The intuition behind the former observation is discussed in the previous experiments (Figure \ref{fig: binomial T reversed}), while the latter observation can be explained as follows. The number of random samples in the truncated data set is relatively small. Hence, a large variance may result in biased estimates of the nominal expected costs that, in turn, results in a higher nominal relative loss for the truncated DRO$_1$ and DRO$_2$ approaches. 


\subsubsection{Results for the unweighted knapsack problem}

\looseness-1 In this section we demonstrate that the obtained practical insights remain valid for a class of unweighted knapsack problems. For the sake of brevity, we focus on a comparison of the baseline DRO approach (\ref{distributionally robust predictor})-(\ref{distributionally robust prescriptor}) and Hoeffding bounds (\ref{hoeffding predictor})-(\ref{hoeffding prescriptor}). Also, the nominal distribution and the distribution of $T_{a}$, $a \in \mathcal{A}$, are supposed to be component-wise binomial and uniform, respectively; recall Tables~\ref{tab: data-generating distributions}~and~\ref{tab: choice of T}. As in the previous section, the parameters $p_a$, $a \in \mathcal{A}$, of the binomial distribution are distributed uniformly over the interval $[0,1]$.

\begin{figure}[htp]
	\begin{center}
		\includegraphics[width = 0.5\linewidth]{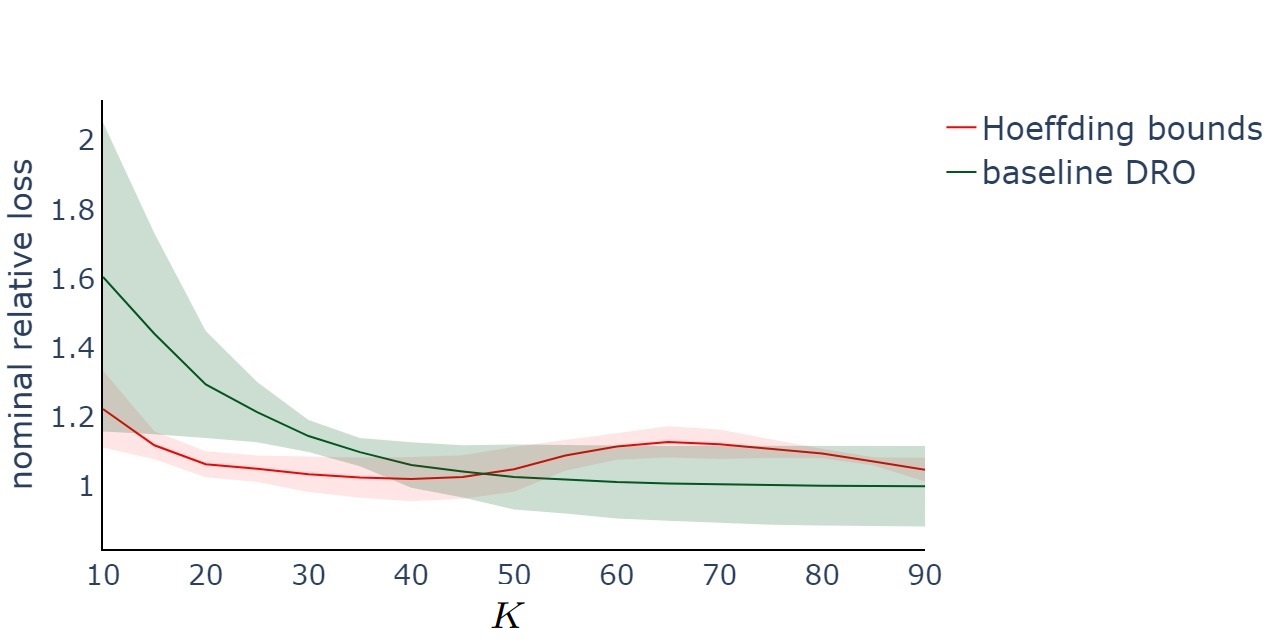}
	\end{center}
	\caption{\scriptsize Average relative loss (\ref{eq: nominal relative loss}) and MADs as a function of $K$, $K \in \{10, 15, \ldots, 90\}$, with $\widetilde{T}_{min} = 10$, $\Delta = 10$ and $n = 100$ under the binomial distribution. The distribution of $T_a$, $a \in \mathcal{A}$, is uniform and $n = 100$. }
	\label{fig: knapsack 1}
\end{figure}

\looseness-1 First, in Figure \ref{fig: knapsack 1} we consider the nominal relative loss (\ref{eq: nominal relative loss}) as a function of $K$ for $\widetilde{T}_{min} = 10$ and $\Delta = 10$; recall definition (\ref{feasible set: knapsack}).
We observe that the baseline DRO approach outperforms Hoeffding bounds only if the parameter $K$ is sufficiently large. 

This observation can be explained as follows. If $K$ is sufficiently small, then we need to select only the components of $\mathbf{c}$ with relatively low expected costs. As outlined in Figure \ref{fig: binomial, costs}, these costs are better estimated using Hoeffding bounds. On the other hand, if $K$ is sufficiently large, than the decision-maker needs to select all components of $\mathbf{c}$ with low expected costs and a part of the components with relatively high expected costs. Since the baseline DRO approach provides better estimates of higher expected costs, we conclude that it performs better for sufficiently large values of $K$. 

\begin{figure}[htp]
	\begin{subfigure}[b]{0.5\textwidth}
		\begin{center}
			\includegraphics[width = \linewidth]{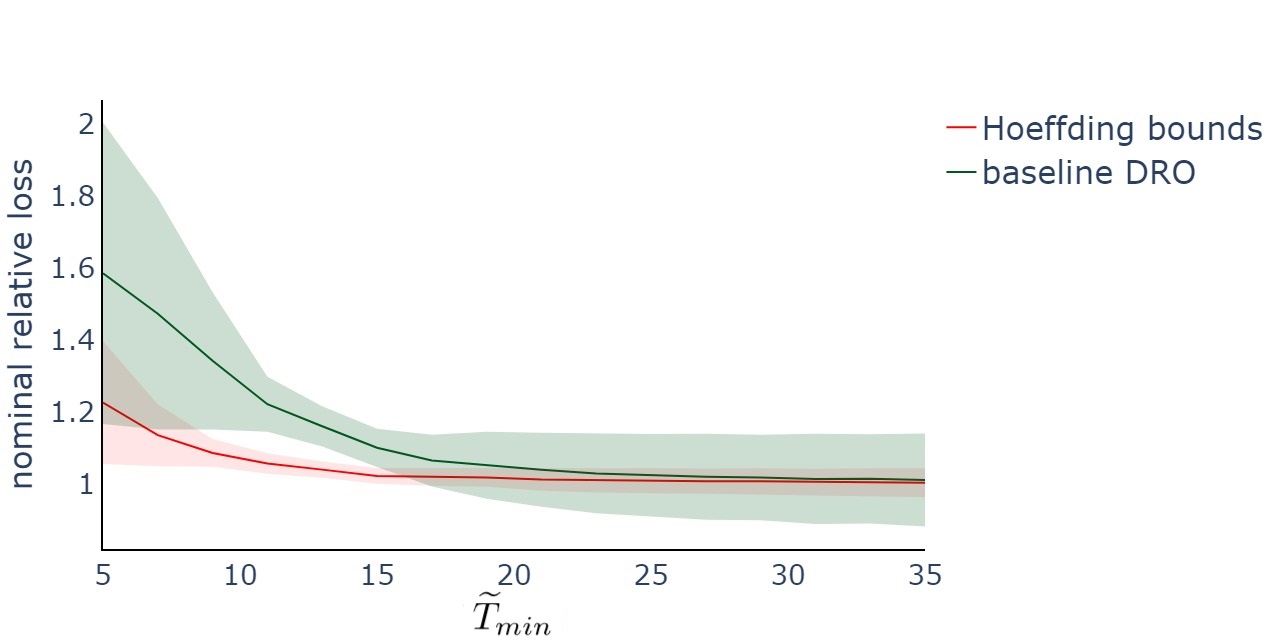}
		\end{center}
		\caption{\scriptsize Binomial distribution and $K = 20$.}
		\label{fig: knapsack 2}
	\end{subfigure}
	\begin{subfigure}[b]{0.5\textwidth}
		\begin{center}
			\includegraphics[width = \linewidth]{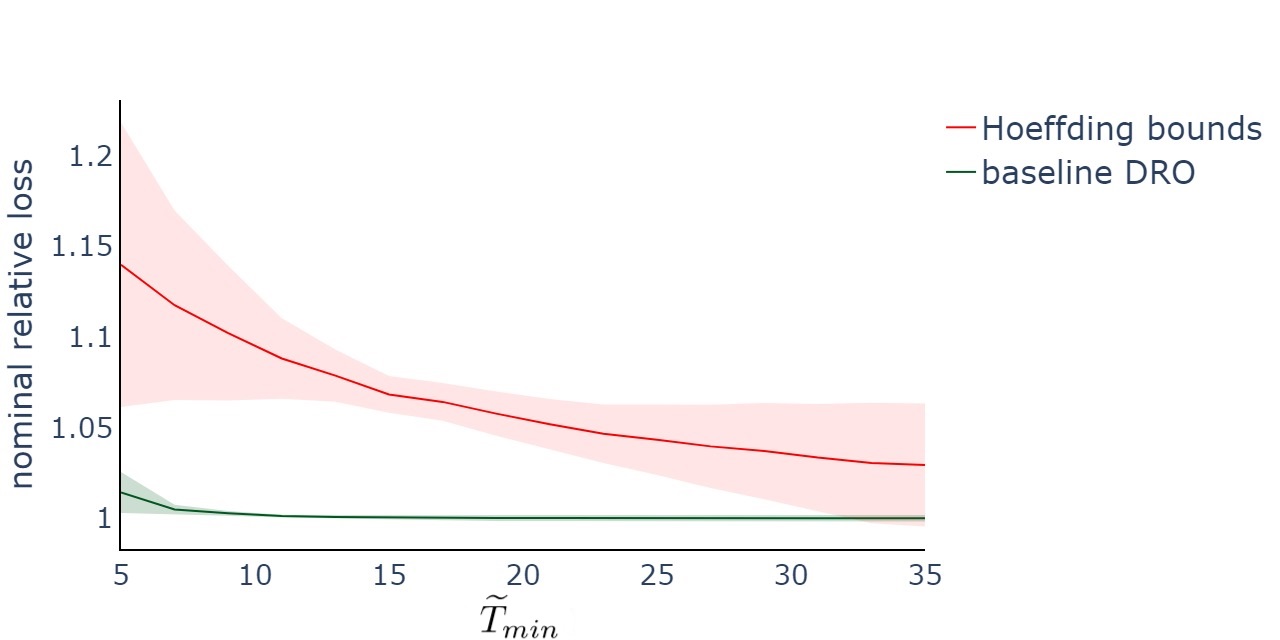}
		\end{center}
		\caption{\scriptsize Binomial distribution and $K = 80$.}
		\label{fig: knapsack 3}
	\end{subfigure}
	\caption{Average relative loss (\ref{eq: nominal relative loss}) and MADs as a function of $\widetilde{T}_{min}$, $\widetilde{T}_{min} \in \{5, 7, \ldots, 35\}$, with $\Delta = 10$, (a) $K = 20$ and (b) $K = 80$ under the binomial distribution. The distribution of $T_a$, $a \in \mathcal{A}$, is uniform and $n = 100$. }
\end{figure}

Finally, we show that the aforementioned results do not depend on the choice of $\widetilde{T}_{min}$. That is, we consider the nominal relative loss (\ref{eq: nominal relative loss}) as a function of $\widetilde{T}_{min}$ with $\Delta = 10$. The plots for $K = 20$ and $K = 80$ are depicted in Figures \ref{fig: knapsack 2} and \ref{fig: knapsack 3}, respectively. 

\looseness-1 \textbf{Summary.} Advantages of the baseline DRO approach (\ref{distributionally robust predictor})-(\ref{distributionally robust prescriptor}) can be summarized as follows. The proposed approach outperforms Hoeffding bounds (\ref{hoeffding predictor})-(\ref{hoeffding prescriptor}) whenever relatively low expected costs can be observed sufficiently often or the gap between the sample sizes for the components of $\mathbf{c}$ is comparatively small. Furthermore, if the variance of the data-generating distribution is sufficiently small, then a truncation of the data set may substantially improve the out-of-sample performance of the baseline DRO approach. In contrast to the optimization model in \cite{VanParys2020}, the outlined effect can be achieved ``for free'', i.e., by leveraging the same DRO model for the truncated data set. 
%

\section{Conclusion} \label{sec: conclusion}
In this paper we consider a class of linear mixed-integer programming problems, where the cost vector is governed by some unknown probability distribution $\mathbb{Q}^*$. The decision-maker attempts to minimize its expected loss under $\mathbb{Q}^*$ using some finite training data set obtained from this distribution. In contrast to the related study of Van Parys~et~al.~\cite{VanParys2020}, we assume that the components of the cost vector are explored not to the same degree. In particular, the proposed modeling approach is motivated by a number of online combinatorial optimization and machine learning problem settings, where the data is collected by trial and errors through multiple decision epochs. 

\looseness-1 For the constructed stochastic programming problem we seek a \textit{prediction rule} that converts the data set into an estimate of the expected value of the objective function and a \textit{prescription rule} that provides an associated estimate of the optimal decision. 
The goal is to find the least conservative prediction and prescription rules, which also hedge against underestimated losses whenever the sample size tends to infinity. We demonstrate that under some mild assumption the associated prediction and prescription problems admit a weakly optimal solution with a number of attractive theoretical properties. First, there is a class of decisions and joint distributions, for which the proposed prediction (prescription) rules are not dominated by any other feasible prediction (prescription) rules. 
Second, the aforementioned solution can be obtained 
by solving a number of component-wise convex distributionally robust optimization problems and a unique instance of the nominal problem. 

We perform numerical experiments, where the out-of-sample performance of the proposed approach is analyzed with respect to several classes of synthetic combinatorial optimization problems. Importantly, we exploit the abovementioned theoretical results to provide some intuition for a numerical validation of our approach. In particular, it turns out that solutions with a reasonably good quality can be obtained whenever the relative difference between sample sizes is sufficiently small or the data set is biased towards the lower expected costs.

Admittedly, our theoretical results exploit a component-wise decomposition of the objective criterion. Therefore, it would be interesting to consider some other risk measures, especially those that do not account for any covariance information among the cost coefficients (but may account for some variance information). To the best of our knowledge, an accurate estimation of the covariance matrix for incomplete data sets is rather challenging; we refer, e.g., to \cite{Liu2019} where this problem is addressed for a particular class of data-generating distributions.
In addition, a rather natural step forward is to extend the current results to the case of a continuous support. We refer the reader to Section 5 in \cite{VanParys2020} for the related discussion in the case of complete data. Finally, it remains an open question whether the proposed approach provides strongly optimal prediction (prescription) rules or not.

\textbf{Acknowledgments.} The authors would like to thank Dr. Oleg Prokopyev for his helpful comments and suggestions. The article was prepared within the framework of the Basic Research Program at the National Research University Higher School of Economics (Sections~ 1-2). The research is funded by RSF project №22-11-00073~(Sections 3-5).

\textbf{Conflict of interest:} The authors declare that they have no conflict of interest.

\bibliographystyle{ieeetr}
\bibliography{bibliography}

\appendix 
\section{Supplementary material} \label{sec: app}
\textbf{The proof of Proposition \ref{proposition 3}.}
Assume that there exists another prediction rule $\hat{f}'$, which is less conservative than $\hat{f}$ and satisfies the asymptotic guarantee (\ref{eq: asymptotic guarantee}). Therefore, for any decision $\mathbf{x} \in X$ and any joint distribution $\mathbb{Q} \in \mathcal{Q}$ there exists $\varepsilon \in \mathbb{R}_{>0}$ such that
\begin{equation} \nonumber
\hat{f}'(\mathbf{x}, \mathbf{Q}) \leq \hat{f}(\mathbf{x}, \mathbf{Q}) - \varepsilon,
\end{equation}
In the remainder of the proof we show that $\hat{f}'$ is infeasible in (\ref{prediction problem}). 

First, we bound the out-of-sample disappointment (\ref{eq: out-of-sample prediction}) from below as:
\begin{equation} \label{eq: WW optimality proof 1} 
\begin{gathered}
\Pr\Big\{f(\mathbf{x}, \mathbb{Q}^*) > \hat{f}'(\mathbf{x}, \widehat{\mathbf{Q}})\Big\} \geq \Pr\Big\{f(\mathbf{x}, \mathbb{Q}^*) > \hat{f}(\mathbf{x}, \widehat{\mathbf{Q}}) - \varepsilon \Big\}, 
\end{gathered}
\end{equation}

Without loss of generality let $z_{a,1} = \overline{z}_a$ for each $a \in \mathcal{A}$.
Since the asymptotic guarantee (\ref{eq: asymptotic guarantee}) must be satisfied for any thinkable data-generating distribution $\mathbb{Q}^* \in \mathcal{Q}$, we may assume that for each $a \in \mathcal{A}$ 
\begin{equation} \nonumber
q^*_{a,i} = 
\begin{cases}
1, \mbox{ if } i = 1,\\
0, \mbox{ otherwise }
\end{cases}\end{equation}

Next, we observe that
\begin{equation} \nonumber
\begin{gathered}
f(\mathbf{x}, \mathbb{Q}^*) = \sum_{a \in \mathcal{A}} \mathbb{E}_{\mathbb{Q}^*_a} \{c_a\} x_a = \sum_{a \in \mathcal{A}} \overline{z}_a x_a > \hat{f}(\mathbf{x}, \mathbf{Q}) - \varepsilon,
\end{gathered}
\end{equation}
for any feasible $\mathbb{Q} \in \mathcal{Q}$ (recall that $\mathbf{Q}$ is a vector of marginal distributions induced by $\mathbb{Q}$),
where the last inequality holds due to the upper bound (\ref{eq: WW optimality upper bound}). Hence,
\begin{equation} \nonumber
\Pr\Big\{f(\mathbf{x}, \mathbb{Q}^*) > \hat{f}'(\mathbf{x}, \widehat{\mathbf{Q}})\Big\} \geq \Pr\Big\{f(\mathbf{x}, \mathbb{Q}^*) > \hat{f}(\mathbf{x}, \widehat{\mathbf{Q}}) -\varepsilon \Big\} = 1
\end{equation} 
and $\hat{f}'$ does not satisfy the asymptotic guarantee (\ref{eq: asymptotic guarantee}). 

\textbf{The proof of Theorem \ref{theorem 3}.}
The proof follows the proof of Theorem \ref{theorem 2} excluding the following minor changes. First, we set $b \in \argmin_{a \in \mathcal{A}}T_a$ and note that the fact $(\widetilde{\mathbf{x}}, \widetilde{\mathbf{Q}}) \in \widetilde{\mathcal{Z}}$ implies that
\begin{equation} \nonumber
\widetilde{q}_{a,i} = \widetilde{q}_{b,i} > 0 \quad \forall i \in \{1, \ldots, d_{min}\}, \; \forall a \in \mathcal{A}, \; a \neq b
\end{equation}
and $\widetilde{q}_{a,i} = 0$ for any $i \in \{d_{min} + 1, \ldots, d_a\}$ and $a \in \mathcal{A}$; recall definition (\ref{eq: set Z 2}). 
Then the expected cost of $c_a$ under the worst-case distribution 
\begin{equation} \nonumber
\mathbb{Q}_a^{(w)} \in \argmax_{\mathbb{Q}_a \in \widehat{\mathcal{Q}}_a}\mathbb{E}_{\mathbb{Q}_a}\{c_a\}
\end{equation}
can be expressed as follows:
\begin{equation}
\begin{gathered}
\mathbb{E}_{\mathbb{Q}^{(w)}_a} \{c_a\} = \sum_{i = 1}^{d_a} q^{(w)}_{a, i}z_{a,i} = \max_{\mathbb{Q}_a \in \widehat{\mathcal{Q}}_a} \sum_{i = 1}^{d_a} q_{a, i}z_{a,i} = \max_{\mathbb{Q}_a \in \widehat{\mathcal{Q}}_a} \sum_{i = 1}^{d_{min}} q_{a, i}z_{a,i} = \\ \max_{\mathbb{Q}_a \in \widehat{\mathcal{Q}}_a} \sum_{i = 1}^{d_{min}} q_{a, i}(w_{b, a} z_{b, i} + v_{b, a}) = w_{b, a} \max_{\mathbb{Q}_a \in \widehat{\mathcal{Q}}_a} \Big( \sum_{i = 1}^{d_{min}} q_{a, i} z_{b, i} \Big) + v_{b, a} \leq \\ w_{b, a} \max_{\mathbb{Q}_b \in \widehat{\mathcal{Q}}_b} \Big( \sum_{i = 1}^{d_{min}} q_{b, i} z_{b, i} \Big) + v_{b, a} = \max_{\mathbb{Q}_b \in \widehat{\mathcal{Q}}_b} \Big( \sum_{i = 1}^{d_b} q_{b, i} z_{a, i} \Big) = \mathbb{E}_{\mathbb{Q}^{(w)}_b} \{c_a\} \nonumber
\end{gathered}
\end{equation}
In contrast to (\ref{eq: weak optimality proof A1'}), we additionally use the definition of the relative entropy ball, i.e., 
\begin{align} \label{eq: app 1}
\sum_{i = 1}^{d_a} q^{(w)}_{a, i}z_{a,i} = & \max \sum_{i = 1}^{d_a} q_{a, i}z_{a,i} \\
\mbox{s.t. } & q_{a, i} > 0 \quad \forall i \in \{1, \ldots, d_a\} \nonumber \\
& \sum_{i = 1}^{d_a} q_{a, i} = 1, \nonumber \\
& \sum_{i = 1}^{d_{min}} \widetilde{q}_{a, i} \ln \frac{\widetilde{q}_{a, i}}{q_{a, i}} \leq r_a   \nonumber
\end{align}
That is, since $\widetilde{q}_{a, i} \ln \frac{\widetilde{q}_{a, i}}{q_{a, i}}$ is a decreasing function of $q_{a, i}$ and $z_{a,1} \geq \ldots \geq z_{a,d_a}$ by construction, the optimal solution of (\ref{eq: app 1}) must satisfy $q^{(w)}_{a,i} = 0$ for $i \in \{d_{min} + 1, \ldots, d_a\}$ and any $a \in \mathcal{A}$. 

Finally, for each $a \in \mathcal{A}$ we define $$\mathbb{Q}^*_a := \lambda \widetilde{\mathbb{Q}}_b + (1 - \lambda) \mathbb{Q}^{(w)}_b$$ for a sufficiently small $\lambda > 0$. In particular, we observe that $q^*_{a,i} = \lambda \widetilde{q}_{a,i} + (1 - \lambda) q^{(w)}_{a,i} > 0$ for any $i \in \{1, \ldots, d_{min}\}$ and $a \in \mathcal{A}$ and $q^*_{a,i} = 0$, otherwise. Thus, in each support set $\mathcal{S}_a$, $a \in \mathcal{A}$, we may leave only the first $d_{min}$ components and then exploit the proof of Theorem \ref{theorem 2} for $d = d_{min}$. This observation concludes the proof.
\end{document}